\documentclass[preprint]{elsarticle}
\usepackage{amssymb}
\usepackage{amsmath}
\usepackage{amstext}
\usepackage{mathptmx}
\usepackage{bm}
\usepackage{booktabs}
\usepackage{multirow, mathtools}
\usepackage{lineno,hyperref}
\usepackage{color}
\usepackage{flafter}
\usepackage{tabularx}
\usepackage[hang]{subfigure}
\usepackage[linesnumbered,lined,ruled]{algorithm2e}

\journal{Journal}

\begin{document}

\begin{frontmatter}

\title{Second order hierarchical partial least squares regression-polynomial chaos expansion for global sensitivity and reliability analyses of high-dimensional models}



\author[mymainaddress]{Ling-Ze Bu\corref{mycorrespondingauthor}}
\cortext[mycorrespondingauthor]{Corresponding author}
\ead{17b933010@stu.hit.edu.cn}
\author[mymainaddress,mysecondaryaddress,mythirdaddress]{Wei Zhao}
\author[mymainaddress,mysecondaryaddress,mythirdaddress]{Wei Wang}

\address[mymainaddress]{School of Civil Engineering, Harbin Institute of Technology, Harbin 150090, China}
\address[mysecondaryaddress]{Key Lab of Structures Dynamic Behavior and Control of the Ministry of Education, Harbin Institute of Technology, Harbin, 150090, China}
\address[mythirdaddress]{Key Lab of Smart Prevention and Mitigation of Civil Engineering Disasters of the Ministry of Industry and Information Technology, Harbin Institute of Technology, Harbin, 150090, China}

\begin{abstract}
	
To tackle the curse of dimensionality and multicollinearity problems of polynomial chaos expansion for analyzing global sensitivity and reliability of models with high stochastic dimensions, this paper proposes a novel non-intrusive algorithm called \textit{second order hierarchical partial least squares regression-polynomial chaos expansion}. The first step of the innovative algorithm is to divide the polynomials into several groups according to their interaction degrees and nonlinearity degrees, which avoids large data sets and reflects the relationship between polynomial chaos expansion and high dimensional model representation. Then a hierarchical regression algorithm based on partial least squares regression is devised for extracting latent variables from each group at different variable levels. The optimal interaction degree and the corresponding nonlinearity degrees are automatically estimated with an improved cross validation scheme. Based on the relationship between variables at two adjacent levels, Sobol' sensitivity indices can be obtained by a simple post-processing of expansion coefficients. Thus, the expansion is greatly simplified through retaining the important inputs, leading to accurate reliability analysis without requirements of additional model evaluations. Finally, finite element models with three different types of structures verified that the proposed method can greatly improve the computational efficiency compared with the ordinary least squares regression-based method.
\end{abstract}


\end{frontmatter}


\section{Introduction}
Mathematical models are essential tools in science and engineering to study and express the behavior of real systems. Exact predictions of system responses can be made only if the model inputs are accurately described. However, this cannot be achieved since uncertainties widely exist. Therefore, uncertainty quantification is crucial for mathematical modeling and has been attracting growing interest around the world.

For real structural systems, loadings, material properties and boundary conditions are usually uncertain and can be described with random variables and random fields. In the uncertainty quantification of structural systems, the analysis of both sensitivity and reliability plays an important role. Sensitivity analysis aims at identifying the effect of different inputs on the uncertainties of model outputs, and includes local and global methods \cite{Blatman2010a}. Global methods are under active research since they can quantitatively describe the relative importance of different inputs by taking their variability in the entire distribution space into account. Various global methods such as variance-based methods \cite{Saltelli2010,Zhang2014}, moment-independent methods \cite{Borgonovo2007,Greegar2015} and derivative-based methods \cite{Sobol2009,Sudret2015} have been proposed in the literature. Reliability analysis aims at evaluating the safety of a structure with failure probability which is the integral of the joint probability density function of model inputs in the failure region. Many reliability analysis methods have been proposed in the past decades including Monte Carlo Simulation (MCS), Importance Sampling \cite{Au1999, Dai2012, Dai2016}, Subset Simulation \cite{Papaioannou2015,Zuev2015}, Line Sampling \cite{Koutsourelakis2004,Koutsourelakis2004a} and other metamodelling approaches such as low rank tensor approximation \cite{Konakli2016}, Neural Network \cite{Dai2015}, Support Vector Machine \cite{Dai2017} and so on.  

The necessary foundation for analyzing sensitivity and reliability of a structural model is uncertainty propagation, i.e. representing the random response. The random response of the model can be explicitly expressed by projecting onto a Hilbert space spanned by the polynomials which are mutually orthogonal with some probability measure. This method is named as \textit{polynomial chaos expansion} \cite{Wiener1938} with three merits:  a sound mathematical background of probability theory and functional analysis, a wide capability suitable for all the random responses with finite variance and an exponential convergence rate for smooth input-output relationships. Currently, it is becoming a popular tool for the uncertainty quantification of structural systems. The probability distribution of stochastic response is characterized by a set of expansion coefficients (i.e. coordinates in the Hilbert space) which can be computed intrusively or non-intrusively. In the intrusive methods, such as the spectral stochastic finite element method (SSFEM) \cite{Ghanem2003}, the coefficients are obtained with a Galerkin scheme which requires modifications of the original deterministic computer codes, hence the term \textit{intrusive} for this approach. On the contrary, the non-intrusive methods only need repeated calls of the deterministic model, which is convenient for analyzing large complex structural systems with well-validated computer codes. Two approaches are usually distinguished, namely projection and regression-based ones. The projection-based approaches are based on the orthogonality between different polynomials, and are essentially to compute multi-dimensional numerical integrals since closed-form solutions are hardly achieved. The regression-based approaches treat the computation of coordinates as a statistical regression problem, having more flexibility in design of experiments and faster convergence rates \cite{Blatman2009}.

The most commonly used regression approach is the ordinary least squares regression (OLSR). However, the required number of model evaluations dramatically increases with the number of model inputs. This problem is called the \textit{curse of dimensionality}. It is crucial to reduce the number of model evaluations because a single simulation of high-fidelity structural model is often expensive and large data sets easily lead to overflow. Moreover, the OLSR will unavoidably encounter \textit{multicollinearity} under small sample sizes, leading to poor estimations of coordinates. A good way to circumvent these problems is to find an advanced regression technique that can be more robust under small sample sizes. Since the randomness of response is mainly caused by a small part of polynomials for most industrial problems, many adaptive sparse PCE approaches based on stepwise regression \cite{Blatman2010,Blatman2010a, Abraham2017, Liu2018}, least angle regression \cite{Blatman2011a}, support vector regression \cite{Cheng2018} D-MORPH regression \cite{Cheng2018a}, compressive sensing \cite{Jakeman2017} were proposed to suffice the curse of dimensionality. 

This paper proposes a novel regression-based PCE methodology for analyzing global sensitivity and reliability of high-dimensional models. Breaking the existing routine that important terms are directly selected from the huge candidate set of finite order polynomials, the new method on one hand builds a hierarchical regression algorithm based on the idea of divide-and-conquer, and on the other hand detects principal directions which captures the probabilistic content in each subspace by using a state-of-the-art regression approach named partial least squares regression (PLSR) \cite{Rosipal2006,Rosipal2010}. Optimal nonlinearity degrees and interaction degrees are automatically selected with a modified cross validation scheme. The new method not only can achieve significant dimension reduction, but contribute to uncovering the latent hierarchical low-dimensional structure of the model as well. This work is distinguished from our previous work \cite{Zhao2019} in two features: Firstly, the manner of subspace division is more elaborate in order to separate the contributions of terms with different degrees of nonlinearity and interaction. Secondly, the updated subspace division manner naturally derives a new regression scheme with an additional level of hierarchy. The optimal subspaces are extracted in a consistent manner in all levels rather than selected with a penalized regression scheme in higher levels.

The remainder of this paper is organized as follows. The next section provides a brief review of polynomial chaos representation of the stochastic response. A novel methodology based on partial least squares regression and hierarchical modeling is proposed in Section 3. The computational gain of the proposed method is illustrated in Section 4 with three finite element models. 

\section{Polynomial chaos representation of the stochastic response}
Denote $(\Omega,\mathcal{F},\mathbb{P}) $ as a probability space and $L^2(\Omega,\mathcal{F},\mathbb{P})$ as the space of random variables with finite second moments:

\begin{equation}
E[X^2]=\displaystyle\int_{\Omega} X^2(\omega) \text{d}\mathbb{P}(\omega)<+\infty
\end{equation}
Assume that the stochastic response of the model $Y(\omega)$ has finite variance, i.e. $Y(\omega) \in L^2(\Omega,\mathcal{F},\mathbb{P})$, thus, it can be represented with the following \textit{polynomial chaos expansion} \cite{Ghanem2003}:

\begin{equation}
Y(\omega)=\overline{Y}+\sum\limits_{p_n=1}^{+\infty} \sum\limits_{\alpha_1+\cdots+\alpha_M=p_n} \beta_{\alpha_1,\dotsc,\alpha_M} \Psi_{\alpha_1,\dotsc,\alpha_M}(\bm{\xi}(\omega)) 
\end{equation}
where $\overline{Y}$ is the mean value, $\bm{\xi}(\omega)=(\xi_1(\omega), \xi_2(\omega),\dotsc,\xi_M(\omega))^{\text{T}}$ is a vector composed of $M$ independent standard Gaussian random variables, each $\Psi_{\alpha_1,\dotsc,\alpha_M}$ is a multidimensional Hermite polynomial with degree $(\alpha_1,\dotsc,\alpha_M)$ and each $ \beta_{\alpha_1,\dotsc,\alpha_M}$ is a deterministic coefficient. For the sake of simplicity,  it is assumed that the random inputs of structures can be described with a set of finite independent random variables. In this context, multidimensional Hermite polynomials can be expressed as tensor product of the corresponding univariate Hermite polynomials. For the case correlated inputs, readers may refer to \cite{Soize2004,Blatman2010,Rahman2017}. The series is mean-square convergent \cite{Cameron1947}, hence can be truncated with maximum order $p_{\mathrm{max}}$:

\begin{equation}\label{maxp}
Y(\omega)=\overline{Y}+\sum\limits_{p_n=1}^{p_{\mathrm{max}}} \sum\limits_{\alpha_1+\cdots+\alpha_M=p_n} \beta_{\alpha_1,\dotsc,\alpha_M} \Psi_{\alpha_1,\dotsc,\alpha_M}(\bm{\xi}(\omega)) 
\end{equation}
Equation (\ref{maxp}) can be written in a compact form:
\begin{equation}
Y(\omega)=\overline{Y}+\sum\limits_{i=1}^{P} \beta_{i} \Psi_{i}(\bm{\xi}(\omega)) 
\end{equation}
where
\begin{equation}\label{P}
P+1=\dfrac{(M+p_{\mathrm{max}})!}{M!p_{\mathrm{max}}!} 
\end{equation}

The key to PCE is the determination of expansion coefficients $\beta_i$ in which the whole information about the probability distribution of the stochastic response is encapsulated. Generally, the computational procedures can be classified into two categories, namely intrusive and non-intrusive methods. In the intrusive methods, the deterministic computer codes have to be modified, which is often cumbersome for large complex structural systems. On the contrary, the deterministic model is treated as a black-box in the non-intrusive method, hence PCE is regarded as a metamodel which can be built with either projection-based or regression-based approaches. The former are based on the orthogonality of polynomials:

\begin{equation}
\beta_i^*=\dfrac{E[\Psi_i Y]}{E[\Psi_i^2]}, i=1,2,\dotsc, P
\end{equation}
 where $E[\cdot]$ is the mean value operator. The denominator can be computed analytically while the nominator need to be estimated with repeated calls of the deterministic model. Computational burden of the latter rapidly increases with the number of input variables, even if Smolyak algorithm \cite{Smolyak1963} can be utilized to alleviate the curse of dimensionality problem. The regression-based approaches are becoming popular due to their flexibility in the design of experiments and fast convergence. Assume $N$ samples are employed to train the metamodel, the vector of coefficients can be computed with equation (\ref{OLS}):
 
\begin{equation}\label{OLS}
\bm{\beta}^*=\bm{\Psi}^+\bm{Y}
\end{equation}
where $\bm{Y}=(Y^{(1)},Y^{(2)},\dotsc,Y^{(N)})^\text{T} $ is the centered response vector, $\bm{\Psi }$ is the polynomial matrix as equation (\ref{Psi}) 
\begin{equation}\label{Psi}
\bm{\Psi}=
\begin{bmatrix}
\Psi_1(\bm{\xi^{(1)}})& \Psi_2(\bm{\xi^{(1)}}) &\cdots & \Psi_{P}(\bm{\xi^{(1)}})\\
\Psi_1(\bm{\xi^{(2)}})& \Psi_2(\bm{\xi^{(2)}}) & \cdots & \Psi_{P}(\bm{\xi^{(2)}})\\
\vdots & \vdots & \vdots & \vdots\\
\Psi_1(\bm{\xi^{(N)}})& \Psi_2(\bm{\xi^{(N)}}) & \cdots & \Psi_{P}(\bm{\xi^{(N)}})\\
\end{bmatrix}
\end{equation}
and $\bm{\Psi}^+$ is the Moore-Penrose generalized inverse of $\bm{\Psi}$. Equation (\ref{OLS}) is the well-known solution of OLSR. However, this approach is effective only if $N>P$. According to equation (\ref{P}), the required number of samples rapidly increases with the stochastic dimension $M$, which is the so-called curse of dimensionality problem. On the other hand, severe multicollinearity will arise under small sample sizes, which leads to poor accuracy. Therefore it is necessary to find an advanced regression technique that can be more robust under small sample sizes. 

\section{Second order hierarchical partial least squares regression-polynomial chaos expansion}\label{prometh}
This section presents a novel hierarchical regression algorithm for constructing PCE based on the PLSR technique which has been applied to deal with multicollinearity in high-dimensional data sets in other science disciplines such as chemometrics \cite{Mehmood2012} and bioinformatics \cite{Cao2008}. According to the High Dimensional Model Representation (HDMR) theory \cite{Rabitz1999b} and the relationship between HDMR and PCE \cite{Sudret2008}, equation (\ref{maxp}) can be rewritten as a summation of component functions with different interaction degrees (i.e. number of inputs in a component function) as equation (\ref{rewritten}):

\begin{equation}\label{rewritten}
\begin{split}
Y(\omega) & =\overline{Y} + \sum\limits_{i=1}^M \sum\limits_{\bm{\alpha}\in \mathcal{I}_{i}}\beta_{\bm{\alpha}}\Psi_{\bm{\alpha}}(\xi_{i}(\omega)) + \sum\limits_{1\leqslant i_1 < i_2 \leqslant M} \sum\limits_{\bm{\alpha}\in \mathcal{I}_{i_1,i_2}}\beta_{\bm{\alpha}}\Psi_{\bm{\alpha}}(\xi_{i_1}(\omega), \xi_{i_2}(\omega)) + \cdots \\
&\quad + \sum\limits_{1\leqslant i_1 <\cdots < i_s \leqslant M} \sum\limits_{\bm{\alpha}\in \mathcal{I}_{i_1,\dotsc,i_s}}\beta_{\bm{\alpha}}\Psi_{\bm{\alpha}}(\xi_{i_1}(\omega), \dotsc,\xi_{i_s}(\omega)) + \cdots \\
&\quad+ \sum\limits_{\bm{\alpha}\in \mathcal{I}_{1,\dotsc,M}}\beta_{\bm{\alpha}}\Psi_{\bm{\alpha}}(\xi_{1}(\omega),\dotsc,\xi_{M}(\omega)) 
\end{split}
\end{equation}
where
\begin{equation}
\mathcal{I}_{i_1,\dotsc,i_s}=\{\bm{\alpha}| \bm{\alpha}\in \mathcal{A}, \alpha_k=0 \iff k \notin (i_1,\dotsc, i_s), \forall k=1,\dotsc, M\}
\end{equation}
\begin{equation}
\mathcal{A}=\left\{\bm{\alpha}\left|\sum\limits_{i=1}^{M}\alpha_i\leqslant p_{\mathrm{max}} \right. \right\} 
\end{equation}
and $\bm{\alpha}=(\alpha_1,\dotsc,\alpha_M)$. This expression coincides with that of polynomial dimension decomposition \cite{Rahman2008}. For most practical problems, the optimal interaction degree $M'$ is much lower than the number of model inputs, in other words, only $M'-$ order HDMR expansion is needed to approximate the input-output relationship with acceptable accuracy. The proposed algorithm aims to adaptively determine the optimal interaction degree and the optimal order of polynomials under each order of interaction degree.

\subsection{Construction of polynomial chaos expansion}

\textit{Step 1}: Initialization

First, select $p_{\mathrm{max}}$ with \textit{prior} knowledge about the nonlinear intensity of the model. Then generate a set of $N$ training samples with the Sobol quasi-random sampling scheme. Next transform the sample points from unit hypercube to the parameter space with isoprobabilistic transform if necessary and run high-fidelity structural model on each sample point to get the corresponding centered model output, denoted as $\bm{F}$. Meanwhile, transform the sample points from unit hypercube to standard Gaussian space and compute the centered polynomial matrix $\bm{E}$.

\textit{Step 2}: Partition of polynomial matrix

The dimension of $\bm{E}$ increases dramatically with the number of model inputs, leading to high computational burden for the OLSR-PCE method. To deal with this issue, we divide $\bm{E}$ into several groups according to their interaction degrees and nonlinearity degrees, and treat them separately and hierarchically. First, $\bm{E}$ is partitioned into the first order subblocks numbered as $p_{\mathrm{max}}$ as equation (\ref{firsub}).

\begin{equation}\label{firsub}
\bm{E}=(\bm{E}^{[1]},\dotsc,\bm{E}^{[p_{\mathrm{max}}]}) 
\end{equation}	
Since $p_{\mathrm{max}} \ll M$ for high dimensional problems,  $M' \leqslant p_{\mathrm{max}}$. $\bm{E}^{[i]} \ (i=1,\dotsc,p_{\mathrm{max}})$ is composed of polynomials which explicitly contains $i$ variables only. Then each $\bm{E}^{[i]}$ is partitioned into the second order subblocks numbered as $(p_{\mathrm{max}}-i+1)$ as equation (\ref{secsub}), 

\begin{equation}\label{secsub}
\bm{E}^{[i]}=(\bm{E}^{[i]}_{[1]},\dotsc,\bm{E}^{[i]}_{[p_{\mathrm{max}}-i+1]}), i=1,\dotsc, p_{\mathrm{max}}
\end{equation}
where $\bm{E}^{[i]}_{[j]} (j=1,\dotsc, p_{\mathrm{max}}-i+1)$ contains polynomials with order $j$ only.

\textit{Step 3}: Hierarchical regression

To determine the optimal order of polynomials which compose the first order component functions of HDMR, let the initial response $\bm{F}_{\text{ini}} =\bm{F}$ and the first order hierarchical partial least squares regression (FOHPLSR) is performed between $\bm{F}_{\text{ini}}$ and $\bm{E}^{[1]}$, as illustrated in Figure \ref{FOHPLSR}:

\begin{figure}[htbp]
	\centering
	\includegraphics[width = 0.7\textwidth]{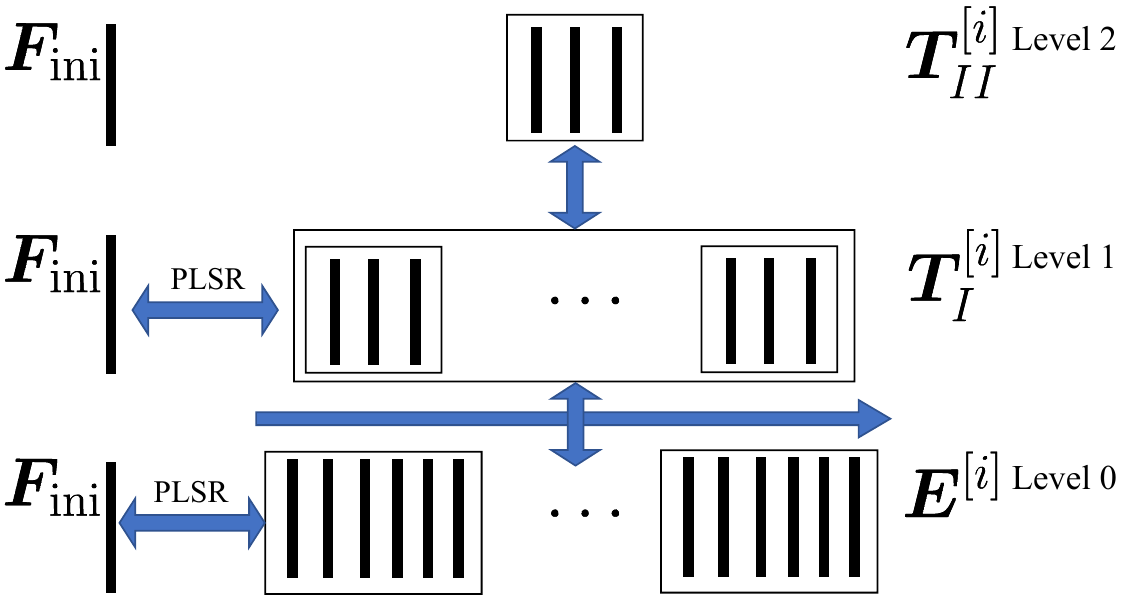}
	\caption{Illustration of first order hierarchical partial least squares regression algorithm}\label{FOHPLSR}
\end{figure}
Let $i=1$, $j=1$, the current response residual and predictor subblock as $\bm{F}_{\text{ini}}$ and $\bm{E}^{[i]}_{[j]}$, respectively. The latent variables of each second order subblock are extracted by performing PLSR between the current response residual and predictor subblock. The fundamental assumption of PLSR is that the behavior of the model is actually controlled by only \textit{a few} factors called \textit{latent variables}. First, latent variables are extracted from $\bm{E}^{[i]}_{[j]}$ and $\bm{F}_{\text{ini}}$, respectively. The latent variables should not only contain the information about variability of the local data sets as much as possible, but also be most correlated. Then, by building regression models among the latent variables, the relationship between $\bm{E}^{[i]}_{[j]}$ and $\bm{F}_{\text{ini}}$ can be built. The latent variables are extracted by iteratively solving the following optimization problem:

\begin{equation}
\left\{\begin{gathered}
\text{Find}\quad \bm{w},\bm{c}\\
\text{s.t.}\quad \max\limits_{\substack{\lVert\bm{w}\rVert_2=1,\\ \lVert\bm{c}\rVert_2=1}}[\text{Cov}^2(\bm{E}^{[i]}_{[j]}\bm{w},\bm{F}_{\text{ini}}\bm{c})]
\end{gathered}\right.
\end{equation}
where $\bm{w}$ and $\bm{c}$ are called loadings. The explanatory and response latent variables are expressed as $\bm{t}=\bm{E}^{[i]}_{[j]}\bm{w}$ and $\bm{u}=\bm{F}_{\text{ini}}\bm{c}$, respectively. The optimization problem can be solved with singular value decomposition algorithm. It is assumed the relationship between $\bm{t}$ and $\bm{u}$ are linear:

\begin{equation}
\bm{u}=b\bm{t}+\bm{H}
\end{equation}
Let $\hat{\bm{u}}=b\bm{t}$, then the contribution of $\bm{t}$ and $\bm{u}$ are deflated from $\bm{E}^{[i]}_{[j]}$ and $\bm{F}_{\text{ini}}$, respectively:
\begin{equation}
\bm{E}^{[i]}_{[j]}=\bm{E}^{[i]}_{[j]}-\bm{tp}^\text{T}
\end{equation}
\begin{equation}
\bm{F}_{\text{ini}}=\bm{F}_{\text{ini}}-\hat{\bm{u}}\bm{q}^\text{T}
\end{equation}
After $h$ iterations, the regression model can be built by using the relationship as equation (\ref{T-E})

\begin{equation}\label{T-E}
\bm{T}=\bm{E}^{[i]}_{[j]}\bm{W}(\bm{P}^\text{T}\bm{W})^{-1}
\end{equation}
where $\bm{T}=(\bm{t}_1,\dotsc, \bm{t}_h)$,  $\bm{W}=(\bm{w}_1,\dotsc, \bm{w}_h)$ and $\bm{P}=(\bm{p}_1,\dotsc, \bm{p}_h)$. In most literature of PLSR, model selection (i.e. selection of the optimal number of iteration) is performed by using the following leave-one-out cross validation error:

\begin{equation}
Err_{\text{LOO}} \equiv \dfrac{1}{N}\sum\limits_{i=1}^{N} \left( F_{\text{ini}}(\bm{\xi}^{(i)}) - \hat{F}_{\text{ini}-i}(\bm{\xi}^{(i)}) \right)^2
\end{equation}
When $Err_{\text{LOO}}$ is lower than a prescribed threshold, such as $1-0.95^2$, the iteration is terminated. However, for PCE, this approach is not only very time-consuming, but also has difficulties in controlling the precision of the model. For OLSR, it can be proven that 

\begin{equation}
F_{\text{ini}}(\bm{\xi}^{(i)}) - \hat{F}_{\text{ini}-i}(\bm{\xi}^{(i)}) = \dfrac{ F_{\text{ini}}(\bm{\xi}^{(i)}) - \hat{F}_{\text{ini}}(\bm{\xi}^{(i)})}{1-h_i}
\end{equation}
where $h_i$ is the $i$th diagonal element of matrix $\bm{E}^{[i]}_{[j]}(\bm{E}_{[j]}^{[i]\text{T}} \bm{E}^{[i]}_{[j]})^{-1} \bm{E}_{[j]}^{[i]\text{T}}$. Similarly, for PLSR, let $h_i$ be the $i$th diagonal element of matrix $ \bm{T}(\bm{T}^{\text{T}} \bm{T})^{-1} \bm{T}^{\text{T}} $, the pseudo leave-one-out cross validation error $Err_{\text{LOO(P)}}$ can be obtained. This is an approximate expression since vector $\bm{w}$ will change when one sample is removed from the original training set, leading to changes of matrix $\bm{T}$. $Err_{\text{LOO(P)}}$ can be normalized by dividing the sample variance, denoted as $\varepsilon_{\text{LOO(P)}}$. However, this error may be unconventional \cite{Blatman2009}, hence the following pseudo cross validation error is introduced for model selection:

\begin{equation}
\varepsilon_{\text{LOO(P)}}^*=\varepsilon_{\text{LOO(P)}} \times \left(1-\dfrac{h}{N}\right)^{-1}\left(1+\text{tr}((\bm{T}^{\text{T}}\bm{T})^{-1})\right)
\end{equation}
The iteration is terminated when $\varepsilon_{\text{LOO(P)}}^*$ reaches its minimum and the corresponding pseudo cross validation error is recorded as $\varepsilon_{\text{LOO(P)},[1]}^{*[1]}$.

Latent variables of other second order subblocks are sequentially obtained by performing PLSR between the current response residual and the predictor subblock. In the process above, latent variables are extracted from level 0 and obtained at level 1. To obtain the latent variables which capture the probabilistic characteristics of univariate polynomials, PLSR is performed sequentially between $\bm{F}_{\text{ini}}$ and $\{\bm{T}^{[1]}_{I,[j]}\}_{j=1}^{k}\ (k\geqslant 2)$, deriving the latent variables at level 2. Meanwhile, the corresponding pseudo cross validation error is recorded as $\varepsilon_{\text{LOO(P)},[k]}^{*[1]}\ (k=2,\dotsc, p_{\mathrm{max}})$. The optimal number of second order subblocks (i.e. the optimal nonlinearity degree) involved in the regression model is selected as the index to the minimum of $\{\varepsilon_{\text{LOO(P)},[k]}^{*[1]}\}_{k=1}^{p_{\mathrm{max}}}$, denoted as $\varepsilon_{\text{LOO(P)}}^{*[1]}$. The process to generate variables at level 2 from level 1 is called a hierarchical operation, deriving the name of the \textit{first order hierarchical partial least squares regression}. The steps of FOHPLSR are summarized in Algorithm \ref{FOHPLSR-al}.

\begin{algorithm}[htbp]
	\caption{The FOHPLSR algorithm}\label{FOHPLSR-al}
	\KwIn{$i$,$\bm{E}^{[i]}$, $\bm{F}_{\text{ini}}$}	
	Let $k=1$, perform PLSR between $\bm{F}_{\text{ini}}$ and $\bm{E}^{[i]}_{[k]}$, get $\bm{T}^{[i]}_{I,[k]}$, $\bm{F}^{[i]}_{[k]}$ and $\varepsilon_{\text{LOO(P)},[k]}^{*[i]}$\\
	\If{$p_{\mathrm{max}}-i+1 \geqslant 2$}{
		\For{$k=2: p_{\mathrm{max}}-i+1$}{perform PLSR between $\bm{F}^{[i]}_{[k-1]}$ and $\bm{E}^{[i]}_{[k]}$, get the explanatory latent variable matrix $\bm{T}^{[i]}_{I,[k]}$ and the corresponding residual $\bm{F}^{[i]}_{[k]}$.
		}		
		\For{$k=2: p_{\mathrm{max}}-i+1$}{perform PLSR between $\bm{F}_{\text{ini}}$ and $(\bm{T}^{[i]}_{I,[1]},\dotsc,\bm{T}^{[i]}_{I,[k]})$, get $\bm{T}^{[i]}_{II,[k]}$ and $\varepsilon_{\text{LOO(P)},[k]}^{*[i]}$}
	}	
	
	Let $k^*=\arg\min\{\varepsilon_{\text{LOO(P)},[k]}^{*[i]}\}_{k=1}^{p_{\mathrm{max}}-i+1}$, $\varepsilon_{\text{LOO(P)}}^{*[i]}=\varepsilon_{\text{LOO(P)},[k^*]}^{*[i]}$, $\bm{F}^{[i]}_{\text{res}}=\bm{F}^{[i]}_{[k^*]}$ and  $\bm{T}^{[i]}_{II}=\bm{T}^{[i]}_{II,[k^*]}$.
	
	\If {$k^*=1$} {$\bm{T}^{[i]}_{II}=\bm{T}^{[i]}_{I,[1]}$}
	\KwOut{$\{\bm{T}^{[i]}_{I,[k]}\}_{k=1}^{k^*}$, $\bm{T}^{[i]}_{II}$,$\bm{F}^{[i]}_{\text{res}}$, $\varepsilon_{\text{LOO(P)}}^{*[i]}$}
\end{algorithm}

The optimal order of polynomials which compose the $i$th $(i=2,\dotsc,p_{\mathrm{max}})$ order component functions of HDMR is estimated by sequentially performing FOHPLSR on the $i$th first order subblock with $\bm{F}_{\text{ini}}=\bm{F}^{[i-1]}_{\text{res}}$. To determine the optimal interaction degree in the regression model, PLSR is performed at level 2 and the corresponding latent variables are obtained at level 3.  Meanwhile, the corresponding pseudo cross validation error is recorded as $\varepsilon_{\text{LOO(P)}}^{*[i]}\  (i=2,\dotsc, p_{\mathrm{max}})$. The optimal interaction degree is selected as the index to the minimum of $\{\varepsilon_{\text{LOO(P)}}^{*[i]}\}_{i=1}^{p_{\mathrm{max}}}$. Finally, the expansion coefficient vector $\bm{\beta}$ can be computed by using equation (\ref{T-E}) from top to bottom level-by-level. The whole process contains two levels of hierarchical operation, hence the formed method is called as \textit{second order hierarchical partial least squares regression-polynomial chaos expansion} (SOHPLSR-PCE) which is summarized as Algorithm \ref{SOHPLSR-al}. The SOHPLSR-PCE algorithm is capable of \textit{group selection} (i.e. selecting a group of polynomials with similar importance) by extracting latent variables, which is effective for dealing with high-dimensional expansions. Thus, the SOHPLSR-PCE method is superior to most approaches including stepwise regression and least angle regression. Meanwhile, the latent variables can be detected on-the-fly (i.e. $\bm{E}^{[i]}_{[j]}$ can be generated sequentially rather than simultaneously), leading to dramatic savings of computer memory.

\begin{algorithm}[htbp]
	\caption{The SOHPLSR-PCE method}\label{SOHPLSR-al}
	\KwIn{$p_{\mathrm{max}}$, $N$}
	Generate $N$ quasi-random samples, compute the corresponding outputs vector $\bm{F}$ and the polynomial matrix $\bm{E}$\\
	Divide $\bm{E}$ into subblocks with equations (\ref{firsub}) and (\ref{secsub})\\
	Let $i=1$, perform FOHPLSR between $\bm{F}$ and $\bm{E}^{[i]}$, get $\bm{T}_{II}^{[i]}$, $\varepsilon_{\text{LOO(P)}}^{*[i]}$ and residual $\bm{F}_{res}^{[i]}$\\
	\If{$p_{\mathrm{max}} \geqslant 2$}{
		\For{$i= 2:p_{\mathrm{max}}$}{perform FOHPLSR between $\bm{F}_{\text{res}}^{[i-1]}$ and $\bm{E}^{[i]}$}
		\For{$i = 2: p_{\mathrm{max}}$}{perform PLSR between $\bm{F}$ and $(\bm{T}_{II}^{[1]},\dotsc,\bm{T}_{II}^{[i]})$, get $\varepsilon_{\text{LOO(P)}}^{*[i]}$ and $\bm{T}^{[i]}_{III}$}		
	}

	Let $i^*=\arg\min\{\varepsilon_{\text{LOO(P)}}^{*[i]}\}_{i=1}^{p_{\mathrm{max}}}$, $\varepsilon_{\text{LOO(P)}}^*=\varepsilon_{\text{LOO(P)}}^{*[i^*]}$, $\bm{F}_{\text{res}}=\bm{F}^{[i^*]}_{\text{res}}$ and  $\bm{T}_{III}=\bm{T}^{[i^*]}_{III}$.\\	
	\If{$i^*=1$}{$\bm{T}_{III}=\bm{T}^{[1]}_{II}$}
	Compute $\bm{\beta}$ by using equation (\ref{T-E}) from top to bottom level-by-level.\\
	\KwOut{$\{\{\bm{T}^{[i]}_{I,[k]}\}_{k=1}^{k^*}\}_{i=1}^{i^*}$, $\{\bm{T}^{[i]}_{II}\}_{i=1}^{i^*}$, $\bm{T}_{III}$, $\varepsilon_{\text{LOO(P)}}^*$, $\bm{\beta}$}
\end{algorithm}

\subsection{Sensitivity analysis}
After obtaining the regression coefficient vector $\bm{\beta}$, similar to the work of \cite{Sudret2008}, variance-based global sensitivity analysis can be performed by a simple post-processing of $\bm{\beta}$, which is expressed in equations (\ref{Si}) and (\ref{STi}):

\begin{equation}\label{Si}
S_{i}= \dfrac{\sum\limits_{\bm{\alpha}\in\mathcal{I}_{i}}\beta_{\alpha}^2E[\Psi_{\alpha}^2(\bm{\xi})]}{\sum\limits_{\bm{\alpha}\in\mathcal{A}}\beta_{\alpha}^2E[\Psi_{\alpha}^2(\bm{\xi})]} 
\end{equation}
\begin{equation}\label{STi}
S_{\text{T}i}=\dfrac{\sum\limits_{\bm{\alpha}\in\mathcal{I}_{i}^*}\beta_{\alpha}^2E[\Psi_{\alpha}^2(\bm{\xi})]}{\sum\limits_{\bm{\alpha}\in\mathcal{A}}\beta_{\alpha}^2E[\Psi_{\alpha}^2(\bm{\xi})]} 
\end{equation}
where $S_i$ and $S_{\text{T}i}$ are called the main and total Sobol indices, respectively.

\subsection{Reliability analysis}
By retaining the input variables whose $S_{\text{T}i}$ values are larger than a prescribed threshold, the PCE metamodel can be reconstructed either with the former samples using OLSR when the cardinality of polynomial set is smaller than the sample size or by using the algorithm introduced above. The failure probability can be estimated with the new metamodel as equation (\ref{PF}),
\begin{equation}\label{PF}
P_{\text{F}}\approx P(g_{\text{PCE}}(\bm{\xi}_{\text{R}})\leqslant 0) = \displaystyle\int_{g_{\text{PCE}}(\bm{x}_{\text{R}})\leqslant 0} f_{\bm{x}_{\text{R}}}(\bm{x}_{\text{R}})\text{d}\bm{x}_{\text{R}}
\end{equation}
where $\bm{\xi}_\text{R}$ is the set of retained input variables, $g_{\text{PCE}}(\bm{x}_{\text{R}})$ is the limit state function of the reconstructed metamodel and $f_{\bm{x}_{\text{R}}}(\bm{x}_\text{R})$ is the joint probability density function of $\bm{\xi}_\text{R}$.

\section{Case Studies}
Three different structures are exampled to examine the accuracy and efficiency of the SOHPLSR-PCE. The finite element models of the structures are firstly built with the MATLAB finite element toolbox \cite{Kattan2008}. Then the computing capacity of the SOHPLSR-PCE is compared with that of the OLSR-PCE. For the latter, it seems that the largest situation is $M=40$ and $p_{\mathrm{max}}=3$ that we can manage with our desk computer. 

\subsection{Simply supported beam}
\begin{figure}[htbp]
	\centering
	\includegraphics[width = 0.6\textwidth]{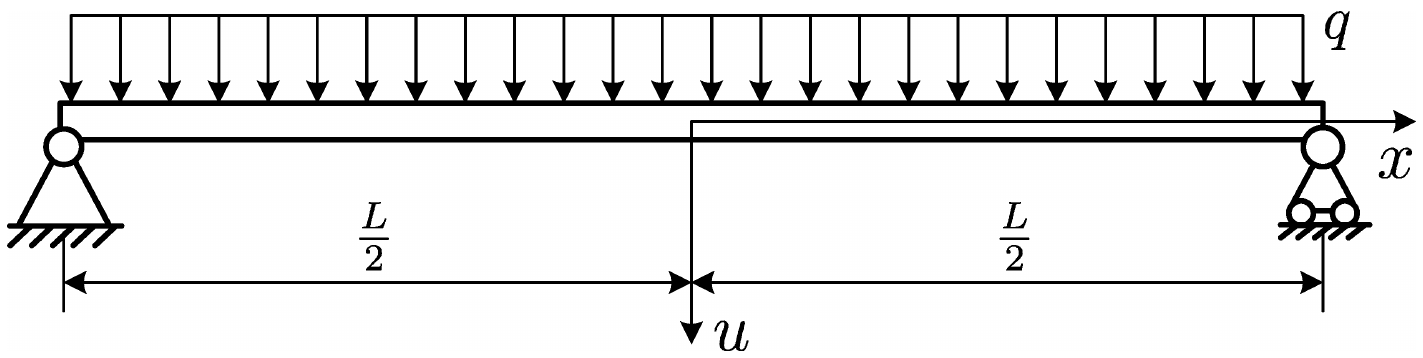}
	\caption{Configuration and loads of a simply supported beam}\label{beam}
\end{figure}

Figure \ref{beam} shows a simply supported beam subjected to an uniformly distributed load.
The beam has length $L=3\text{m}$ and inertial moment $I=8\times10^{-6}\text{m}^4$.  The intensity of distributed loading is $q=13\text{kN/m}$. The elastic modulus of the beam is represented with a non-Gaussian random field as in equation (\ref{Ex1_E})
\begin{equation}\label{Ex1_E}
E(x,\omega)=\exp (N(x,\omega))
\end{equation}
where $N(x,\omega)$ is a homogeneous Gaussian random field whose Pearson product moment correlation coefficient function $\rho  (x,x')$ is 
\begin{equation}
\rho  (x,x')=\exp \left(-\dfrac{|x-x'|}{l}\right) 
\end{equation}
where correlation length $l=0.5\text{m}$. $N(x,\omega)$ is discretized with the first 40 components of Karhunen-Lo\`{e}ve expansion. The mean and coefficient of variation of elastic modulus are $\mu_E=210\text{GPa}$ and $\delta_E=0.2$, respectively. The quantity of interest is the vertical midspan displacement, denoted as $u$. The beam is discretized with 100 elements with equal length. The failure event is defined as $u>$0.012m.

\textit{Step 1}: Construction of polynomial chaos expansion

According to the steps introduced in Section \ref{prometh}, first, the highest order of polynomial $p_{\mathrm{max}}$ is selected as 3, and the corresponding number of polynomial is $P=\text{C}_{40+3}^{3}-1=12340$. Define $\gamma=N/P$ as the expanded sample ratio which is selected as 0.05, 0.2, 0.6, 1 and 2, respectively, and the corresponding sample size $N$ is 617, 2468, 7403, 14808 and 24680, respectively. For each $\gamma$, the metamodel is built with the OLSR-PCE. To illustrate fast convergence of the SOHPLSR-PCE comparing with the OLSR-PCE, in another group of experiments, we define the raw sample ratio $\phi$=2, 4, 6, 8, 10 and 12 ($\gamma$=0.0065, 0.0130, 0.0194, 0.0259, 0.0324 and 0.0389, N=80, 160, 240, 320, 400 and 480), respectively, and built the metamodel with the SOHPLSR-PCE under each $\phi$.

\textit{Step 2}: Global sensitivity analysis

\begin{figure}[htbp]
	\centering
	\noindent\makebox[\textwidth][c]{
	\subfigure[OLSR-PCE]
	{\includegraphics[width=0.5\textwidth]{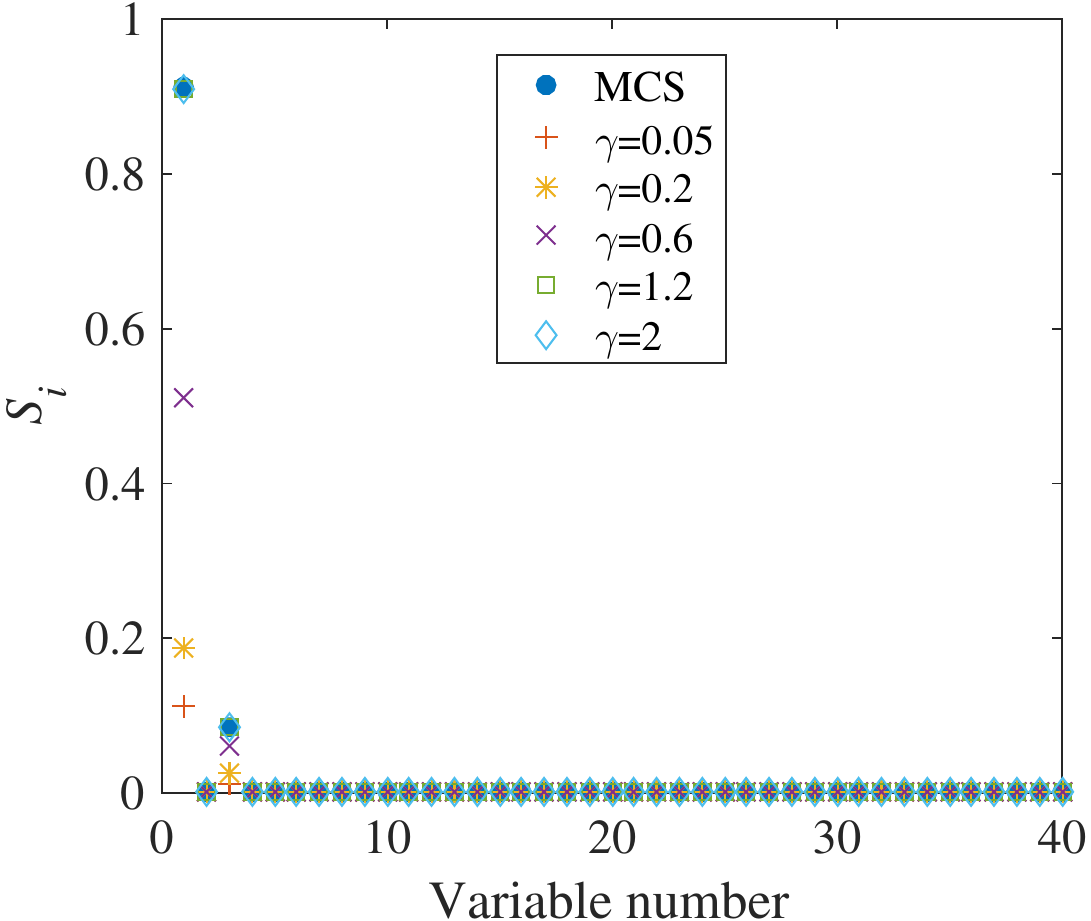}}
	\subfigure[SOHPLSR-PCE]
	{\includegraphics[width=0.5\textwidth]{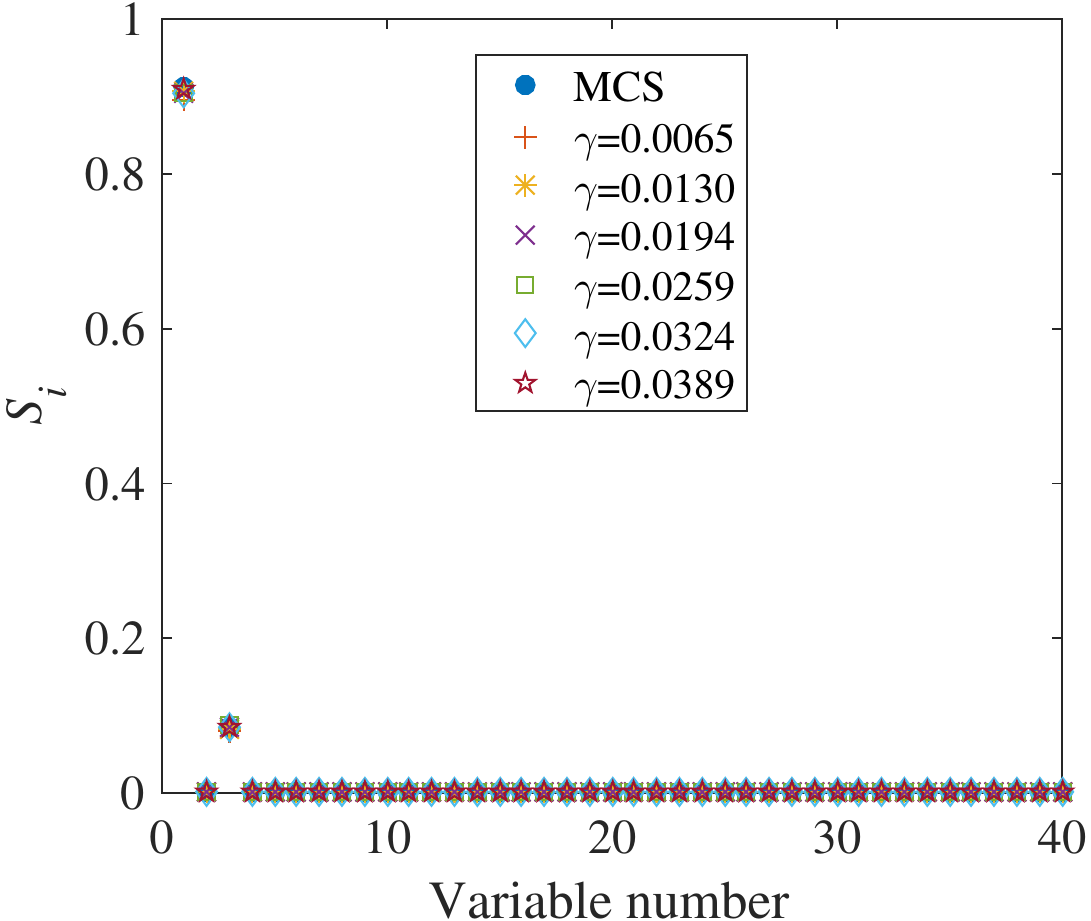}}}
	\caption{Comparison of the main Sobol indices}\label{Ex1_SiOLSHPLS}
	\vspace{0em}
\end{figure}

\begin{figure}[htbp]
	\centering
	\noindent\makebox[\textwidth][c] {
	\subfigure[OLSR-PCE]
	{\includegraphics[width=0.5\textwidth]{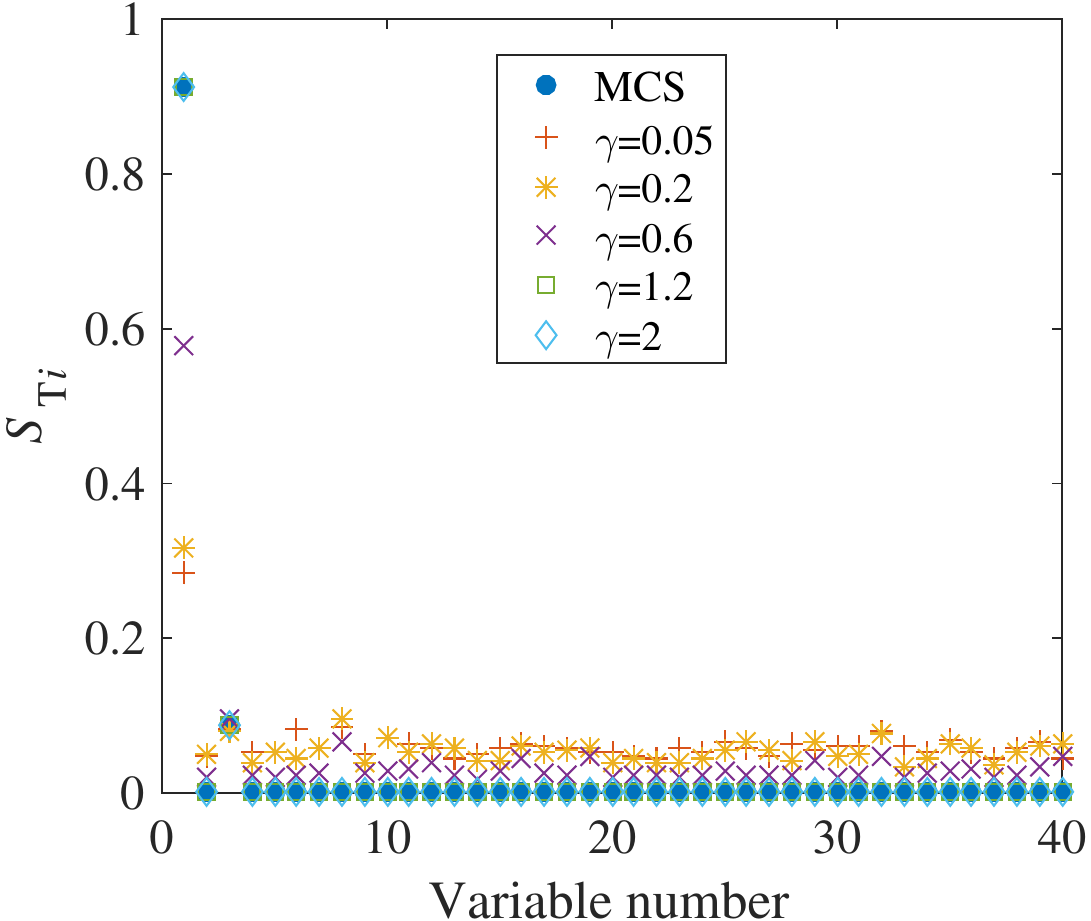}}
	\subfigure[SOHPLSR-PCE]
	{\includegraphics[width=0.5\textwidth]{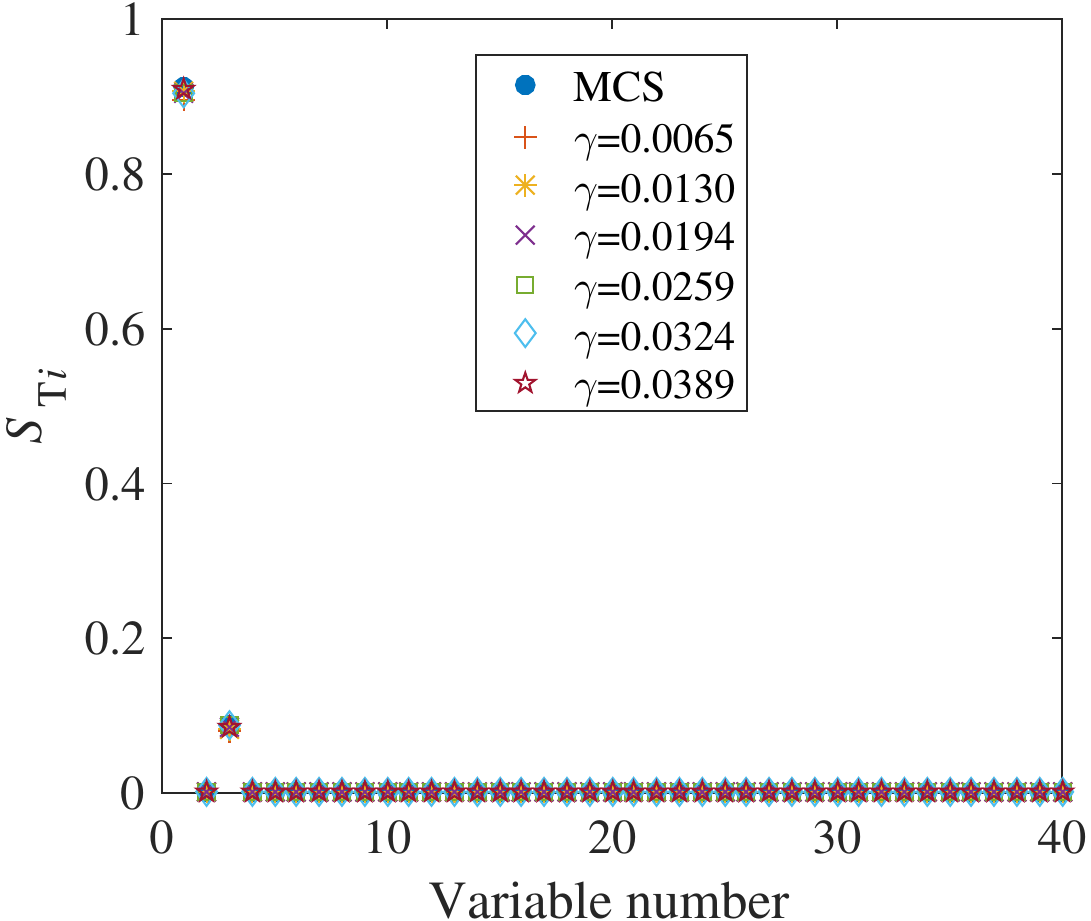}}}
	\caption{Comparison of the total Sobol indices}\label{Ex1_STOLSHPLS}
	\vspace{0em}
\end{figure}

\begin{figure}[htbp]
	\centering
	\noindent\makebox[\textwidth][c]{
	\subfigure[The main Sobol indices]
	{\includegraphics[width=0.5\textwidth]{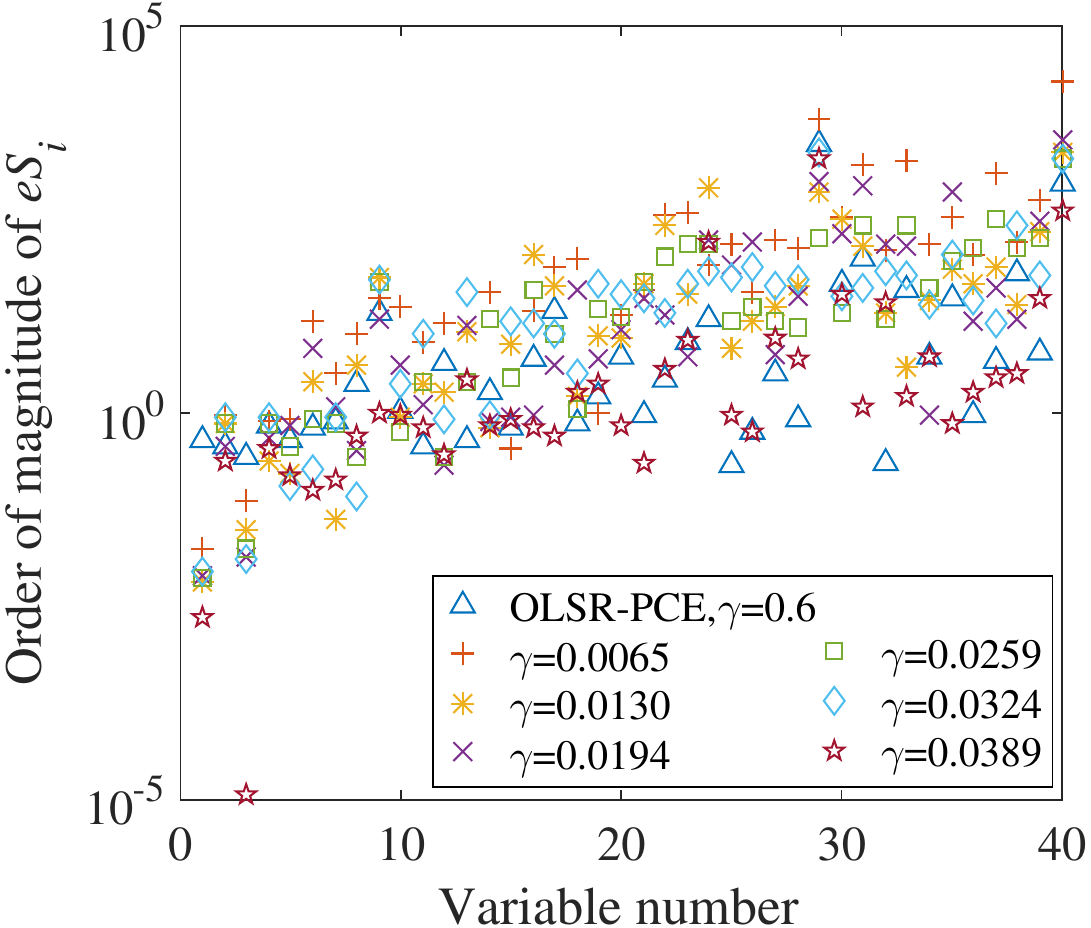}}
	\subfigure[The total Sobol indices]
	{\includegraphics[width=0.5\textwidth]{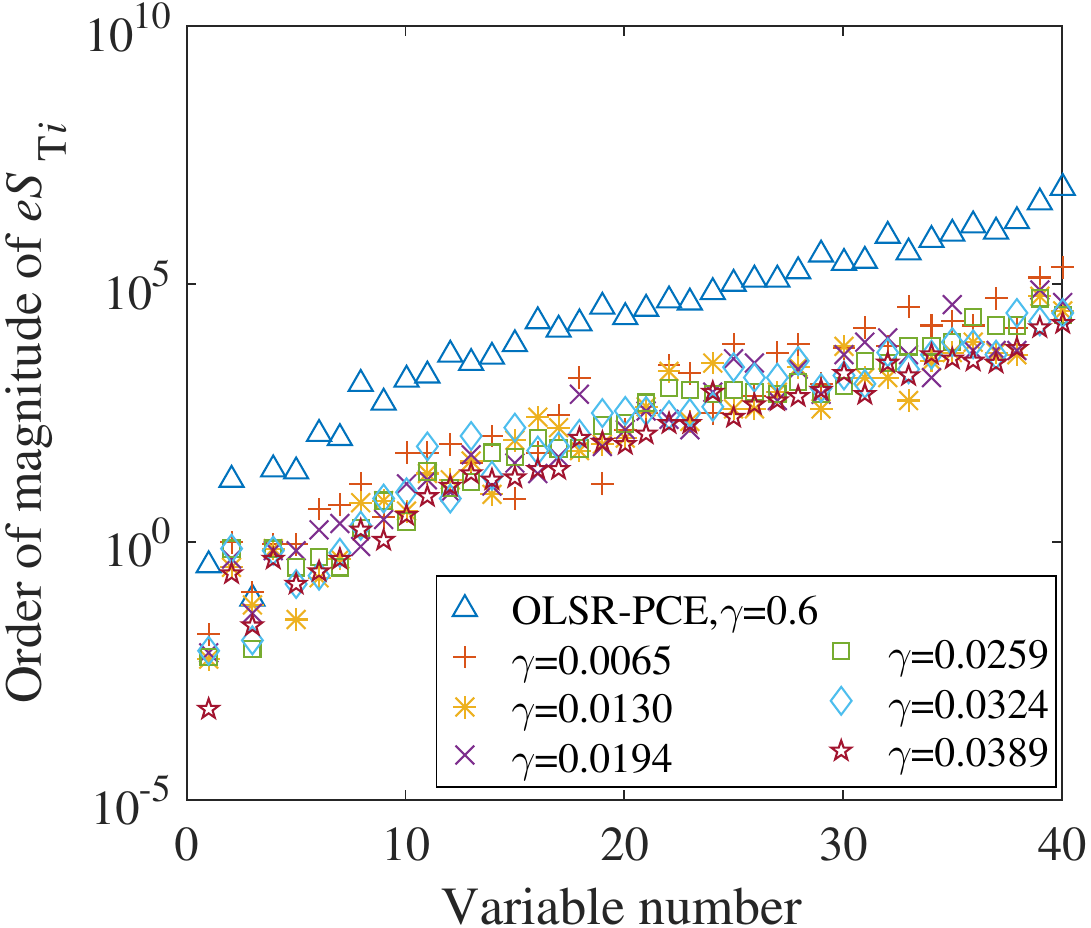}}}
	\caption{Distributions of order of magnitudes of relative errors of Sobol indices}\label{Ex1_eHPLS}
	\vspace{0em}
\end{figure}

The reference solution of each Sobol index is obtained by using MCS with $2\times 10^6$ samples. Comparisons of the main and total Sobol indices are illustrated in Figures \ref{Ex1_SiOLSHPLS} and \ref{Ex1_STOLSHPLS}.
It can be seen from Figures \ref{Ex1_SiOLSHPLS} and \ref{Ex1_STOLSHPLS} that the randomness of variables 1 and 3 have significant impact on the randomness of the model output while the impact of other variables can be ignored. Also, the impact of interaction of different variables is very weak. For the OLSR-PCE, the main Sobol indices of the important variables is inaccurate when $\gamma$ is very low such as 0.05. The solution can qualitatively reflect the distribution of $S_i$ values (i.e. reflect right ranking) only if $\gamma \geqslant 0.2$. To accurately describe the distribution of $S_i$ values, $\gamma$ must be higher than 1.2. The accuracy of the total Sobol indices increases with the sample size slower than that of the main Sobol indices because the former are more dependent on the coefficients of cross terms. Accuracy of the individual coefficients cannot be ensured for OLSR-PCE under small sample ratios, leading to overestimation of the impact of cross terms, which in turn leads to wrong rankings of the importance of the input variables. On the contrary, for both main and total Sobol indices, the SOHPLSR-PCE can consistently get accurate results under all sample ratios, even at very low expanded sample ratio $\gamma=0.0065$, hence the computational efficiency is 185 times as many as the OLSR-PCE. This is because the optimal interaction degree in the expansion and the corresponding nonlinearity degrees are automatically selected by the proposed method, and multicollinearity is dramatically alleviated by using PLSR. For more detailed comparison, we denote the relative error of the main and total Sobol indices as $eS_i$ and $eS_{\text{T}i}$, respectively, and compute the distributions of orders of magnitude for $eS_i$ and $eS_{\text{T}i}$ under each sample ratio using the SOHPLSR-PCE and under expanded sample ratio 0.6 using the OLSR-PCE, respectively, as shown in Figure \ref{Ex1_eHPLS}. 
Because the Sobol indices of some variables are very close to zero, the results computed with the SOHPLSR-PCE will not lead to wrong judgments of important variables and their rankings, although the relative errors of Sobol indices of these variables are at the level of $10^0$-$10^4$ orders of magnitude. For the main Sobol indices, the proposed method does not have significant advantage over the traditional one since $S_i$ values of the latter have already reached high levels of accuracy except for variables 1 and 3. While for the total Sobol indices, by using the SOHPLSR-PCE, the order of magnitudes of relative errors are 1-2 orders of magnitude lower than those of OLSR-PCE whose computational costs are 15-90 times as many as the former. Therefore, for global sensitivity analysis, the SOHPLSR-PCE outperforms the traditional one in terms of computational cost and accuracy.

To uncover the hierarchical structure of the regression model in detail, the numbers of latent variables at different levels are listed in Table \ref{HieStru}.
\begin{table}
	\centering
	\caption{Number of latent variables in each level}\label{HieStru}
	\begin{tabular*}{\textwidth}{c@{\extracolsep{\fill}}c@{\extracolsep{\fill}}c@{\extracolsep{\fill}}c@{\extracolsep{\fill}}c@{\extracolsep{\fill}}c@{\extracolsep{\fill}}c@{\extracolsep{\fill}}c@{\extracolsep{\fill}}c@{\extracolsep{\fill}}c@{\extracolsep{\fill}}c}
		\toprule
		\multirow{2}*{$\gamma$} & \multicolumn{10}{c}{Number of latent variables}\\
		\cmidrule(lr){2-11} & (1,1,1) & (1,1,2) & (1,1,3) & (1,2,2) & (1,2,3) & (1,3,3) & (2,1) & (2,2) &(2,3) & (3)\\
		\midrule
		0.0065&7&	5&	3&	35&	1&	1&	6&	14&	1&	10\\
		0.0130&5&	3&	2&	58&	3&	0&	3&	21&	0&	11\\
		0.0194&4&	3&	2&	67&	1&	0&	3&	16&	0&	8\\
		0.0259&3&	3&	2&	86&	1&	1&	3&	16&	1&	9\\
		0.0324&3&	3&	2&	108&	1&	0&	2&	17&	0&	7\\
		0.0389&4&	3&	1&	138&	1&	0&	2&	19&	0&	7\\
		\bottomrule 	   
	\end{tabular*}
\end{table}
In this table, (1,2,3) means the latent variables extracted from the third order polynomials in the second order HDMR expansion at level 1, (2,1) means the latent variables extracted from the first order HDMR expansion at level 2, (3) means the latent variables extracted at level 3, and so on. It can be seen in Table \ref{HieStru} that the first order HDMR component functions can be represented with a sum of 3-7 linear combinations of first order polynomials, 3-5 linear combinations of second order univariate polynomials and 1-3 linear combinations of third order univariate polynomials. The total number of latent variables at level 1 is 8-15 while the univariate polynomials with orders 1-3 are numbered as 120. Similar comments can be made for the other orders of HDMR component functions. Thus, the regression model is more parsimonious with the increasing of variable level, which fits the fundamental assumption of PLSR. The optimal interaction degree is 2 or 3 and the optimal order of polynomials under each interaction degree is 3. The numbers are not deterministic due to the randomness of each experiment. The optimal numbers of latent variables at different variable levels are nearly the same except for the bivariate second order polynomials. However, it is unimportant since this part has little impact on the variability of model output.

\textit{Step 3}: Reliability analysis

\begin{figure}[htbp]
	\centering
	\includegraphics[width = 0.5\textwidth]{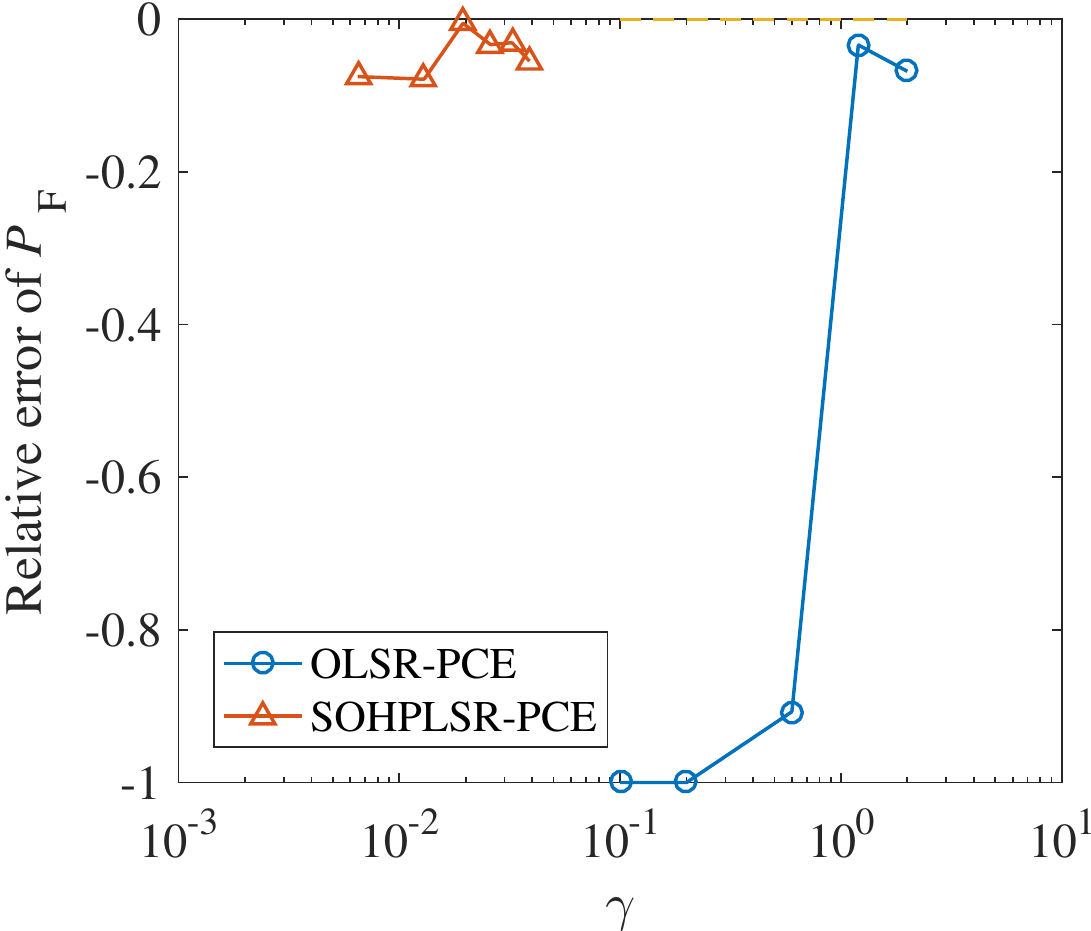}
	\caption{Comparison of relative errors of failure probabilities}\label{Ex1_ePF_HPLS_PEM}
\end{figure}

The threshold for screening the active random inputs is set as 0.018. Then the PCE metamodel is reconstructed with the retained variables based on the results of global sensitivity analysis. The reference result of failure probability is 0.0023, which is obtained by using MCS with $3\times 10^6$  samples. The relative errors under different sample ratios are illustrated in Figure \ref{Ex1_ePF_HPLS_PEM}.
It can be seen in Figure \ref{Ex1_ePF_HPLS_PEM} that the SOHPLSR-PCE can provide rather stable results with relative error less than 10\% when $\gamma \geqslant 0.0065$. Since the proposed method can compute the Sobol indices with high accuracy and low computational cost, the retained variables 1 and 3 can be effectively detected, thus the effective stochastic dimension is reduced to two, and the metamodel can be easily reconstructed with OLSR. On the contrary, the OLSR-based method tends to overestimate the contribution of cross terms under small sample sizes, thus leads to overestimation of total Sobol indices, so that more unimportant variables are retained in the reconstructed metamodel, finally leads to overfitting of the metamodel. Therefore, in this case, $\gamma \geqslant 1.2$ is available to gain acceptable accuracy and stability. Hence, the SOHPLSR-PCE outperforms the OLSR-PCE with a computational gain factor 185.

\subsection{Plane truss}
\begin{figure}[htbp]
	\centering
	\includegraphics[width = 0.7\textwidth]{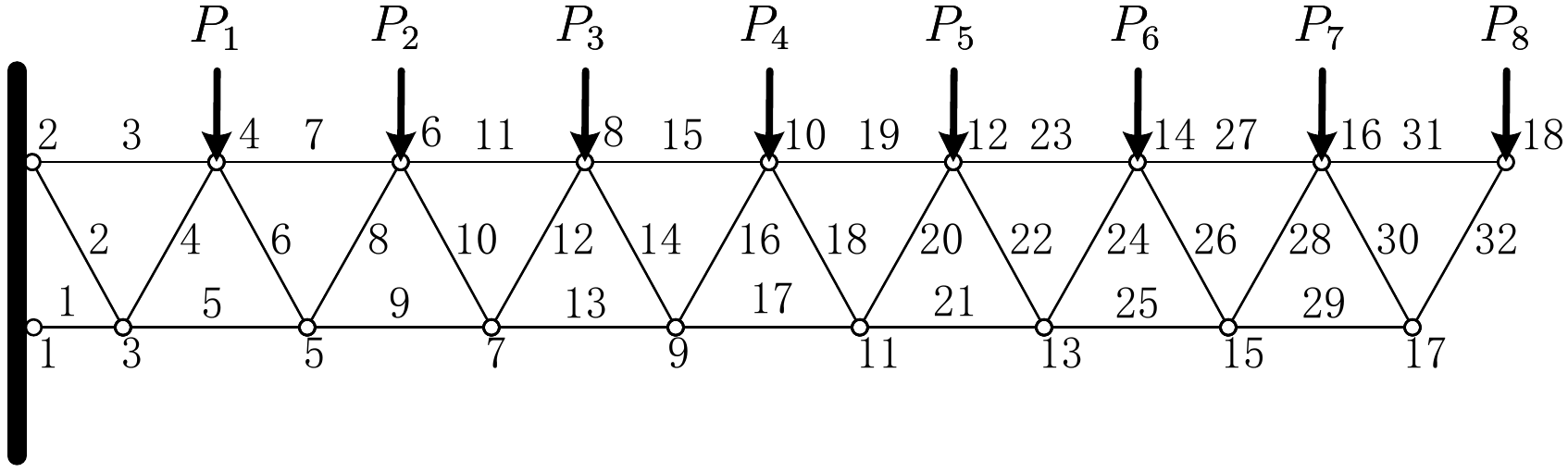}
	\caption{Configuration and loads of a plane truss}\label{truss}
\end{figure}

Figure \ref{truss} shows a plane truss subjected to vertical loads.
Each bar has elastic modulus $E_i\ (i=1,\dotsc,32)$ and diameter 20mm. All the inputs are independent random variables and their distribution parameters are listed in Table \ref{tbeg3}.  The quantity of interest is the vertical displacement of node 18, denoted as $u$. The failure event is defined as $u>$0.210m.
\begin{table}[htbp]
	\caption{Distribution parameters of the inputs}\label{tbeg3}
	\begin{tabular}{cccc}
		\toprule
		Variables & Distribution type & Mean & Standard deviation  \\
		\midrule
		$E_i\ (i=1,\dotsc,32)\ (\text{Pa})$  & Lognormal & $2.0 \times 10^{11} $ & $3.0 \times 10^{10}$\\
		$P_1\ (\text{N})$  & Extreme 1 & $1.2\times 10^4 $ & $2.0 \times 10^3 $\\
		$P_2\ (\text{N})$  & Extreme 1 & $1.0\times 10^4 $ & $1.5 \times 10^3 $\\
		$P_3\ (\text{N})$  & Extreme 1 & $9.0 \times 10^3  $ & $1.2 \times 10^3 $\\
		$P_j\ (j=4,\dotsc,8)\ (\text{N})$  & Extreme 1 & $8.0 \times 10^3 $ & $1.0 \times 10^3 $\\
		\bottomrule
	\end{tabular}
\end{table}

\textit{Step 1}: Construction of polynomial chaos expansion

According to the steps introduced in Section \ref{prometh}, the highest order of polynomial $p_{\mathrm{max}}$ is firstly selected as 3, and the corresponding number of polynomial is $P=\text{C}_{40+3}^{3}-1=12340$.  First, $\gamma$ is selected as 0.05, 0.2, 0.6, 1 and 2 ($N$=617, 2468, 7403, 14808 and 24680), respectively. For each $\gamma$, the OLSR-PCE is used for building the metamodel. To illustrate fast convergence of the SOHPLSR-PCE comparing with the OLSR-PCE, in another group of experiments, $\phi$ is selected as 2,4, 6, 8, 10 and 12 ($\gamma$=0.0065, 0.0130, 0.0194, 0.0259, 0.0324 and 0.0389, $N$=80, 160, 240, 320, 400 and 480), respectively. The metamodel is built with the SOHPLSR-PCE under each $\phi$. 

\textit{Step 2}: Global sensitivity analysis

\begin{figure}[htbp]
	\centering
	\noindent\makebox[\textwidth][c]{	
	\subfigure[OLSR-PCE]
	{\includegraphics[width=0.5\textwidth]{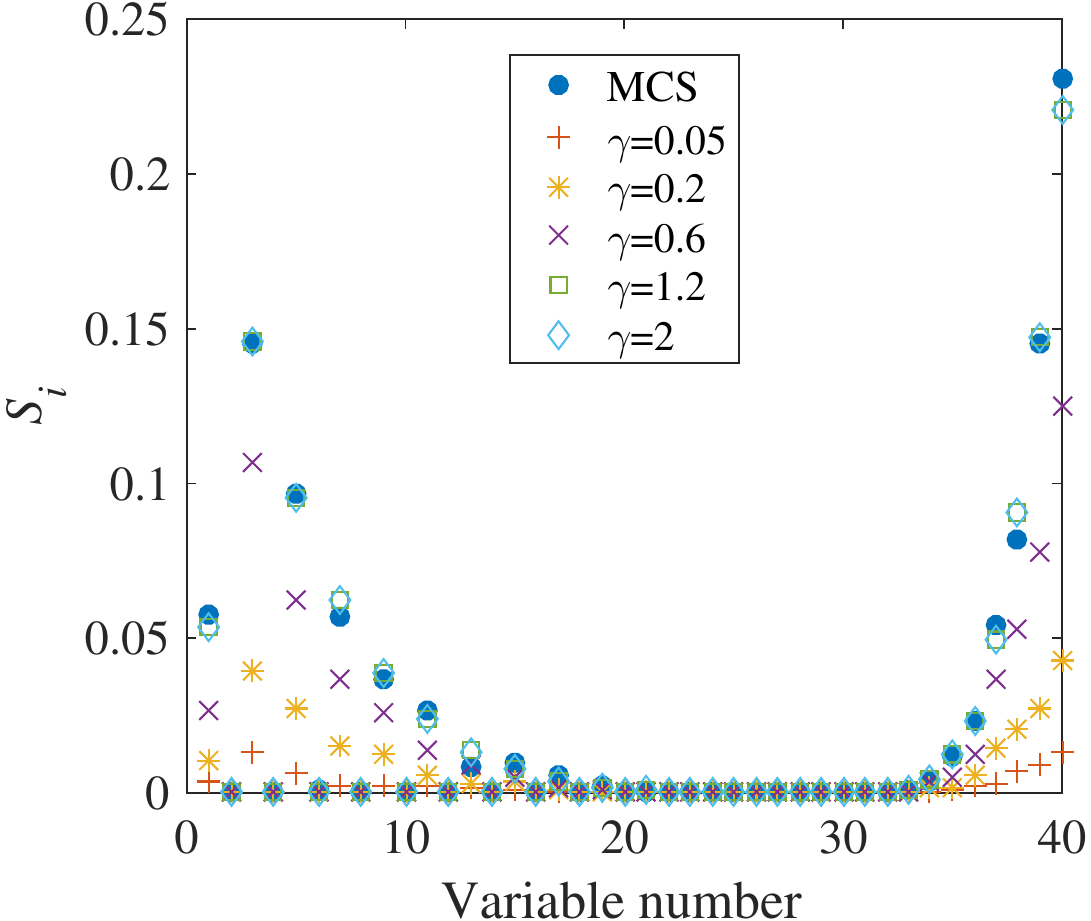}}
	\subfigure[SOHPLSR-PCE]
	{\includegraphics[width=0.5\textwidth]{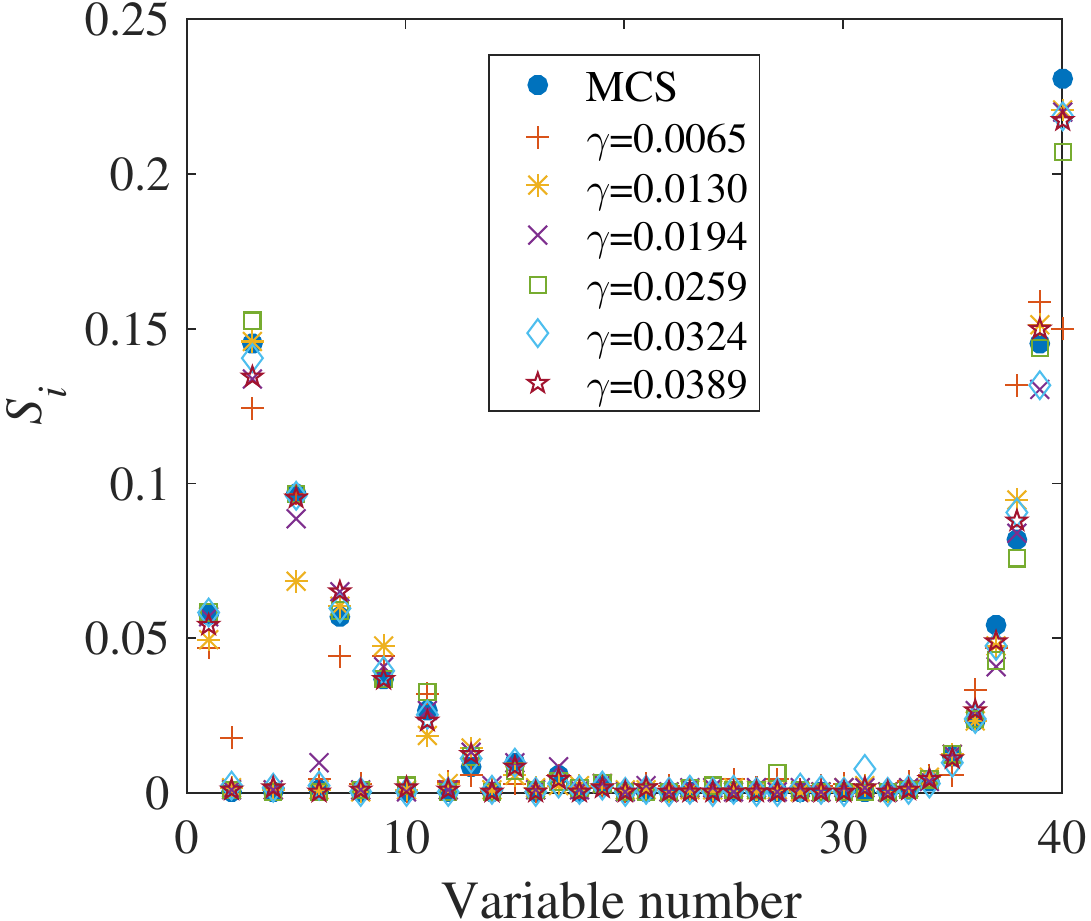}}}
	\caption{Comparison of the main Sobol indices}\label{Ex3_SiOLSHPLS}
	\vspace{0em}
\end{figure}

\begin{figure}[htbp]
	\centering
	\noindent\makebox[\textwidth][c]{	
	\subfigure[OLSR-PCE]
	{\includegraphics[width=0.5\textwidth]{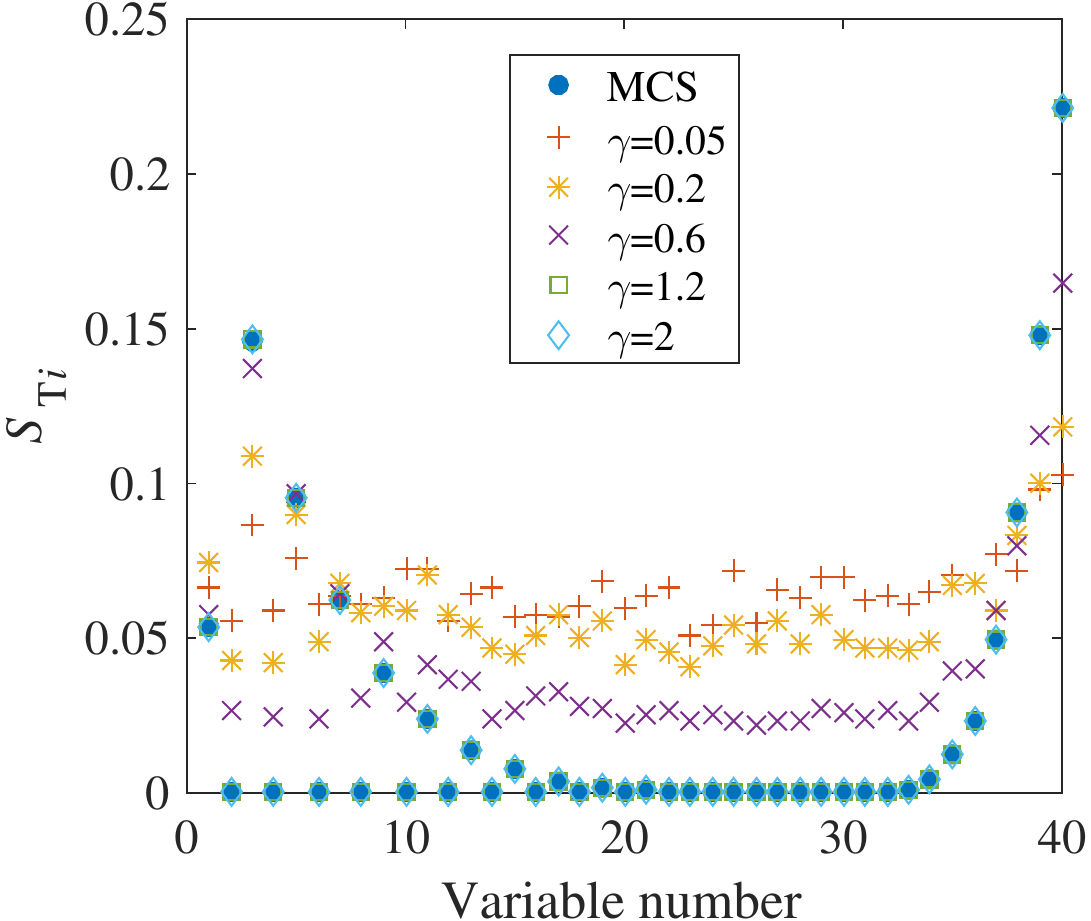}}
	\subfigure[SOHPLSR-PCE]
	{\includegraphics[width=0.5\textwidth]{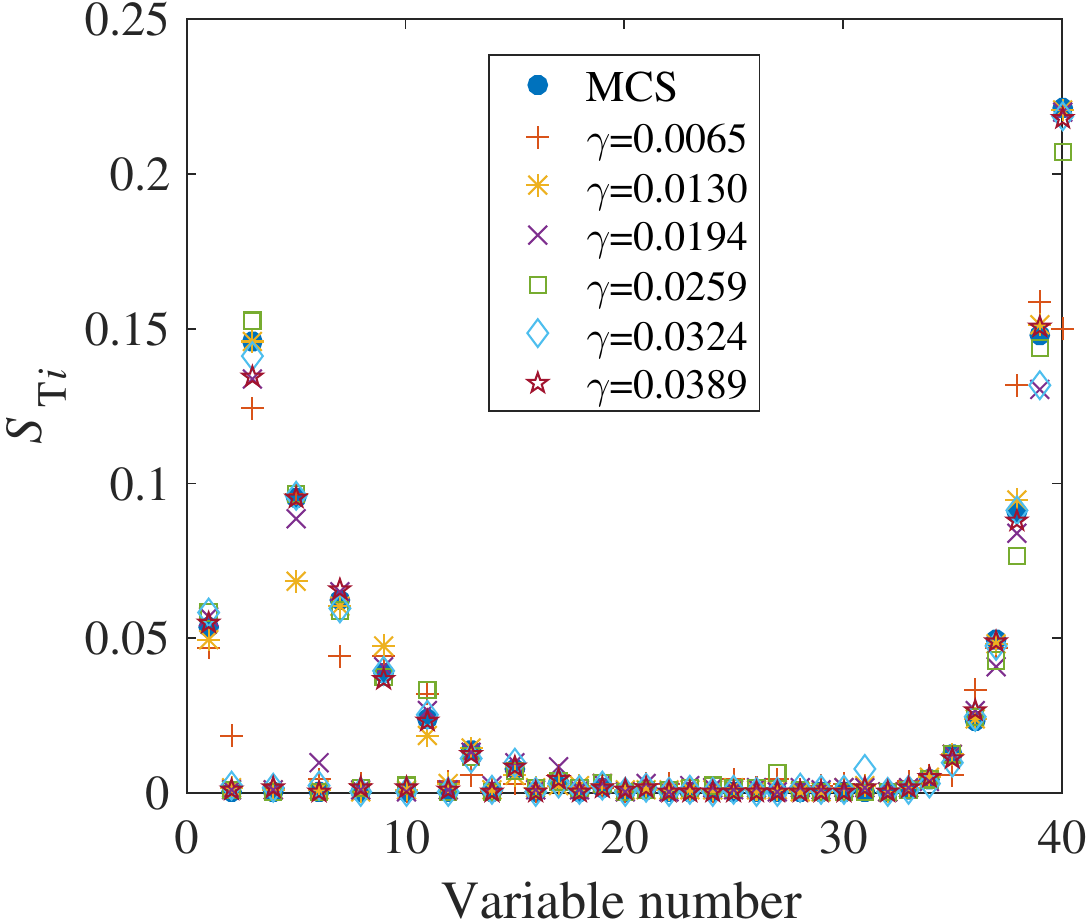}}}
	\caption{Comparison of the total Sobol indices}\label{Ex3_STOLSHPLS}
	\vspace{0em}
\end{figure}

The reference solution of each Sobol index is obtained by using MCS with $2\times 10^6$ samples. Comparisons of the main and total Sobol indices are illustrated in Figures \ref{Ex3_SiOLSHPLS} and \ref{Ex3_STOLSHPLS}.
It can be seen that variables 1, 3, 5, 7, 9, 11, 13, 15, 35-40 individually have significant impact on the uncertainty of model output, but their interactions have a weak impact. Since OLSR-PCE tends to overestimate the contribution of cross terms, the $S_i$ values are generally smaller and the $S_{\text{T}i}$ values are generally greater than the reference solution. The discrepancies become lower with the increasing of the sample ratio. Results with acceptable accuracy can be obtained only when $\gamma \geqslant 1.2$. However, accuracy can be considerably improved by using the SOHPLSR-PCE with noticeably smaller $\gamma$. As the first example, it is acceptable that large errors occur for some variables because their impact is very weak and the errors will not affect the screening and ranking of important variables. As expected, accuracy and stability of the results increase with $\gamma$. Only $\gamma=0.0130$ is needed to quantitatively describe the distribution of the main and total Sobol indices on the whole and only $\gamma=0.0389$ is needed to get highly accurate results. Therefore, the SOHPLSR-PCE outperforms the OLSR-PCE in terms of required number of model evaluations.

\textit{Step 3}: Reliability analysis

\begin{figure}[htbp]
	\centering
	\includegraphics[width = 0.5\textwidth]{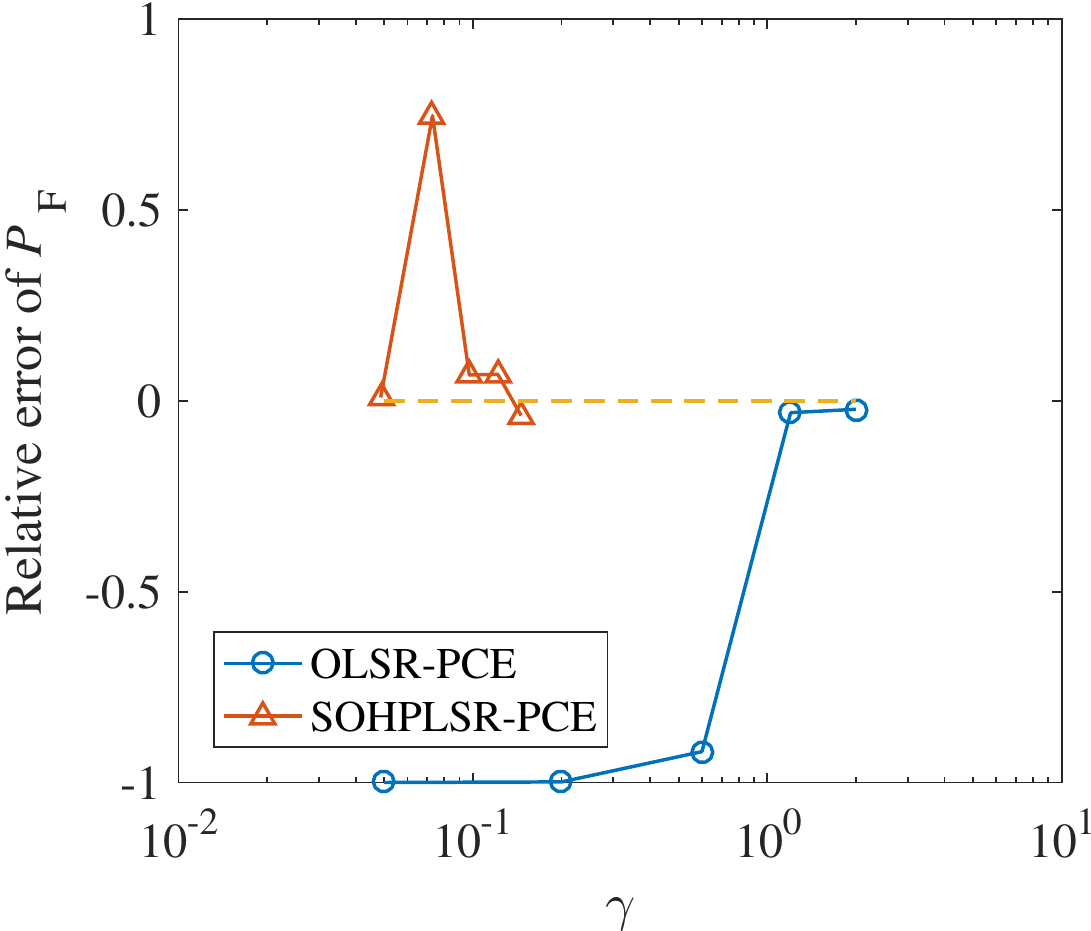}
	\caption{Comparison of relative errors of failure probabilities}\label{Ex3_ePF_HPLS_PEM}
\end{figure}

 The threshold for screening the important random inputs is set as 0.005. Numerical experiments show that the metamodels reconstructed under the original sample ratios cannot provide results with satisfactory stability. To get more accurate and stable results, $\phi$ is increased to 15, 22.5, 30, 37.5 and 45 ($\gamma$=0.0486, 0.0729, 0.0972, 0.1216 and 0.1459, $N$=600, 900, 1200, 1500 and 1800), respectively. Then the failure probabilities are computed with the algorithm introduced in Section \ref{SOHPLSR-al} under the new sample sets. The reference result of failure probability is 0.0052, which is obtained by using MCS with $3\times 10^6$ samples. The relative errors under different sample ratios are illustrated in Figure \ref{Ex3_ePF_HPLS_PEM}.
The global sensitivity analysis above show that the proposed method can exactly detect the important inputs with low computational cost. However, the number of important variables is around 15 which is much more than that in the first example. Therefore it is difficult to ensure the stability and accuracy of the results computed with the metamodel constructed under the original sample ratios, indicating the accuracy and stability of the SOHPLSR-PCE can be improved by enriching the sample size. As illustrated in Figure \ref{Ex3_ePF_HPLS_PEM}, the errors begin to converge to less than 10\% when $\gamma \geqslant 0.0972$. Although the efficiency decreases compared with the first example, it is still 12 times as many as the OLSR-PCE which needs $\gamma\geqslant 1.2$ to keep accuracy.

\subsection{Spatial truss}
\begin{figure}[htbp]
	\centering
	\includegraphics[width = 0.5\textwidth]{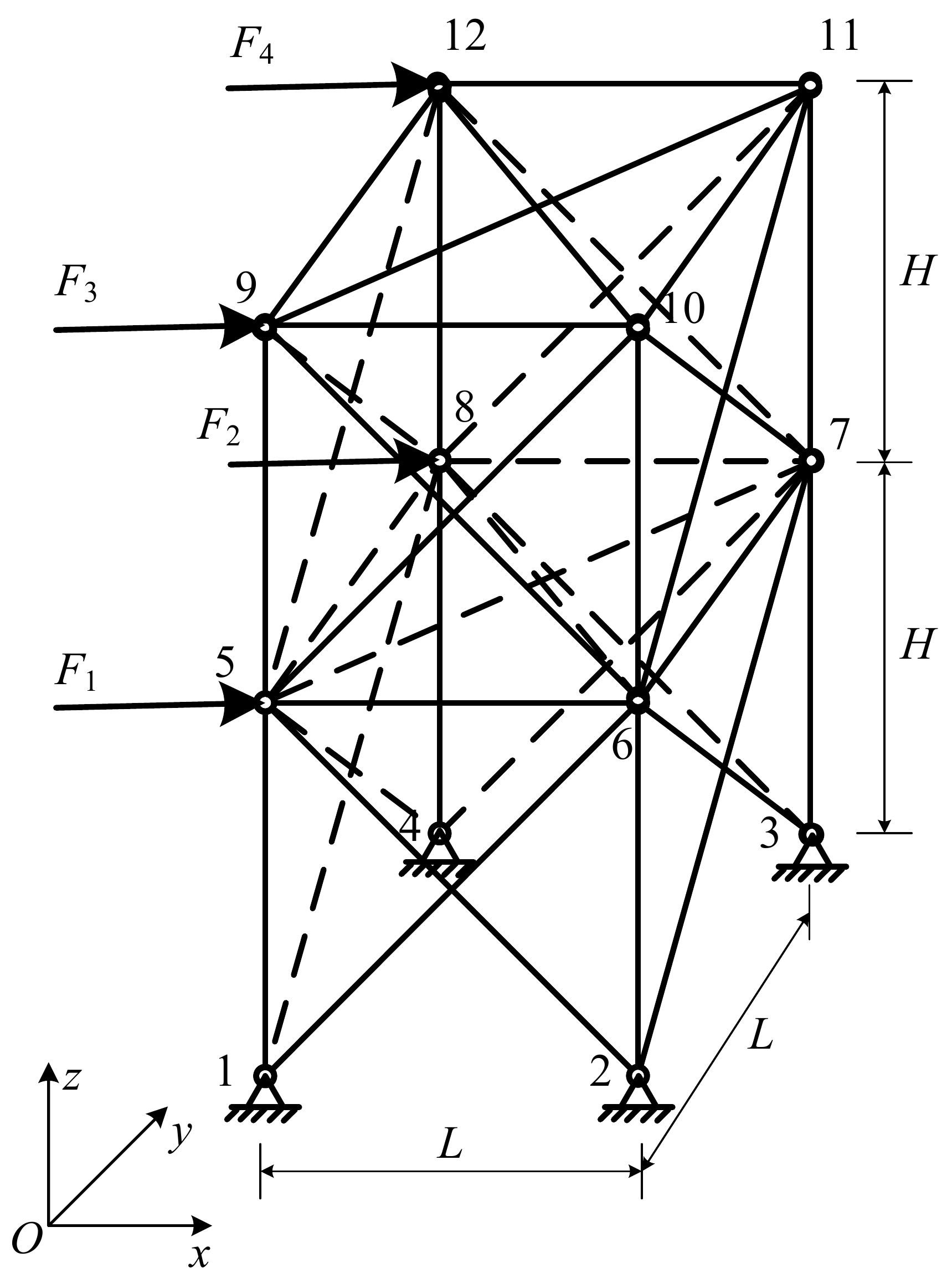}
	\caption{Configuration and loads of a spatial truss}\label{spatialtruss}
\end{figure}

Figure \ref{spatialtruss} shows a spatial truss subjected to horizontal loads.
Each bar (with number listed in Table \ref{elenumspatialtruss}) has elastic modulus $E_i\ (i=1,\dotsc,36)$ and diameter 14mm. $L$ and $H$ are set to be 1m, respectively. All the inputs are independent random variables and their distribution parameters are listed in Table \ref{tbeg6}. The quantity of interest is the maximum horizontal displacement of the top nodes, denoted as $u$. The failure event is defined as $u>$0.004m. 
\begin{table}[htbp]
	\centering
	\caption{Element numbers of the spatial truss}\label{elenumspatialtruss}
	\begin{tabular}{cccccc}
	    \toprule
	    Elements & Node 1 & Node 2 & Elements & Node 1 & Node 2 \\
	    \midrule
	    1 & 1 & 5 & 19 & 5 & 9 \\
	    2 & 5 & 4 & 20 & 9 & 8 \\
	    3 & 5 & 2 & 21 & 9 & 6 \\
	    4 & 2 & 6 & 22 & 6 & 10 \\
	    5 & 6 & 1 & 23 & 10 & 5 \\
	    6 & 6 & 3 & 24 & 10 & 7 \\
	    7 & 3 & 7 & 25 & 7 & 11 \\
	    8 & 7 & 2 & 26 & 11 & 6 \\
	    9 & 7 & 4 & 27 & 11 & 8 \\
	    10 & 4 & 8 & 28 & 8 & 12 \\
	    11 & 8 & 3 & 29 & 12 & 7 \\
	    12 & 8 & 1 & 30 & 12 & 5 \\
	    13 & 5 & 6 & 31 & 9 & 10 \\
	    14 & 6 & 7 & 32 & 10 & 11 \\
	    15 & 7 & 8 & 33 & 11 & 12 \\
	    16 & 8 & 5 & 34 & 12 & 9 \\
	    17 & 5 & 7 & 35 & 9 & 11 \\
	    18 & 6 & 8 & 36 & 10 & 12 \\ 
	    \bottomrule
	\end{tabular}
\end{table}

\begin{table}[htbp]
	\centering
	\caption{Distribution parameters of the inputs}\label{tbeg6}
	\begin{tabular}{cccc}
		\toprule
		Variables & Distribution type & Mean & Standard deviation  \\
		\midrule
		$E_i\ (i=1,\dotsc,36)\ (\text{Pa})$  & Lognormal & $2.0 \times 10^{11} $ & $3.0 \times 10^{10}$\\
		$F_i\ (i=1,\dotsc,4)\ (\text{N})$  & Extreme 1 & $1.0\times 10^4 $ & $1.5\times 10^3 $\\
		\bottomrule
	\end{tabular}
\end{table}

\textit{Step 1}: Construction of polynomial chaos expansion

In the last example, we compare the effects of selecting different values of $p_{\mathrm{max}}$. Values of $p_{\mathrm{max}}$ are selected as 2 and 3, respectively. The corresponding values of $P$ are 12340 and 860. Following the steps introduced in Section \ref{prometh}, in each case, metamodels are built by using OLSR-PCE with $\gamma$= 0.05, 0.2, 0.6, 1.2, 2 ($N$= 43, 172, 516, 1032, 1720 for $p_{\mathrm{max}}$=2, $N$= 617, 2468, 7403, 14808, 24680 for $p_{\mathrm{max}}$=3) and SOHPLSR-PCE with $\phi$= 2, 4, 6, 8, 10, 12 ($N$= 80, 160, 240, 320, 400, 480; $\gamma$= 0.093, 0.186, 0.2791, 0.3721, 0.4651, 0.5581 for $p_{\mathrm{max}}$=2, $\gamma$=0.0065, 0.0130, 0.0194, 0.0259, 0.0324, 0.0389 for $p_{\mathrm{max}}$=3), respectively. 

\textit{Step 2}: Global sensitivity analysis

\begin{figure}[htbp]
	\centering
	\noindent\makebox[\textwidth][c]{	
		\subfigure[OLSR-PCE]
		{\includegraphics[width=0.5\textwidth]{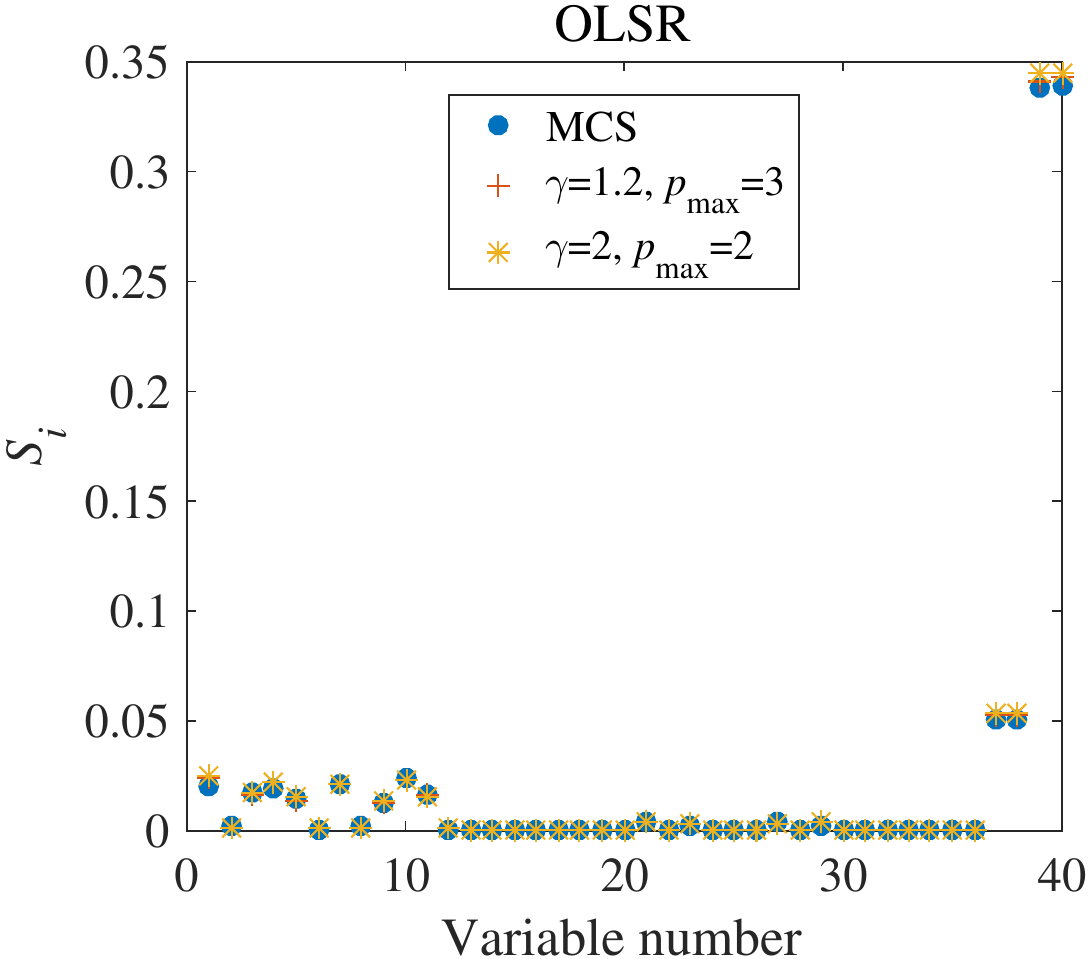}}
		\subfigure[SOHPLSR-PCE]
		{\includegraphics[width=0.5\textwidth]{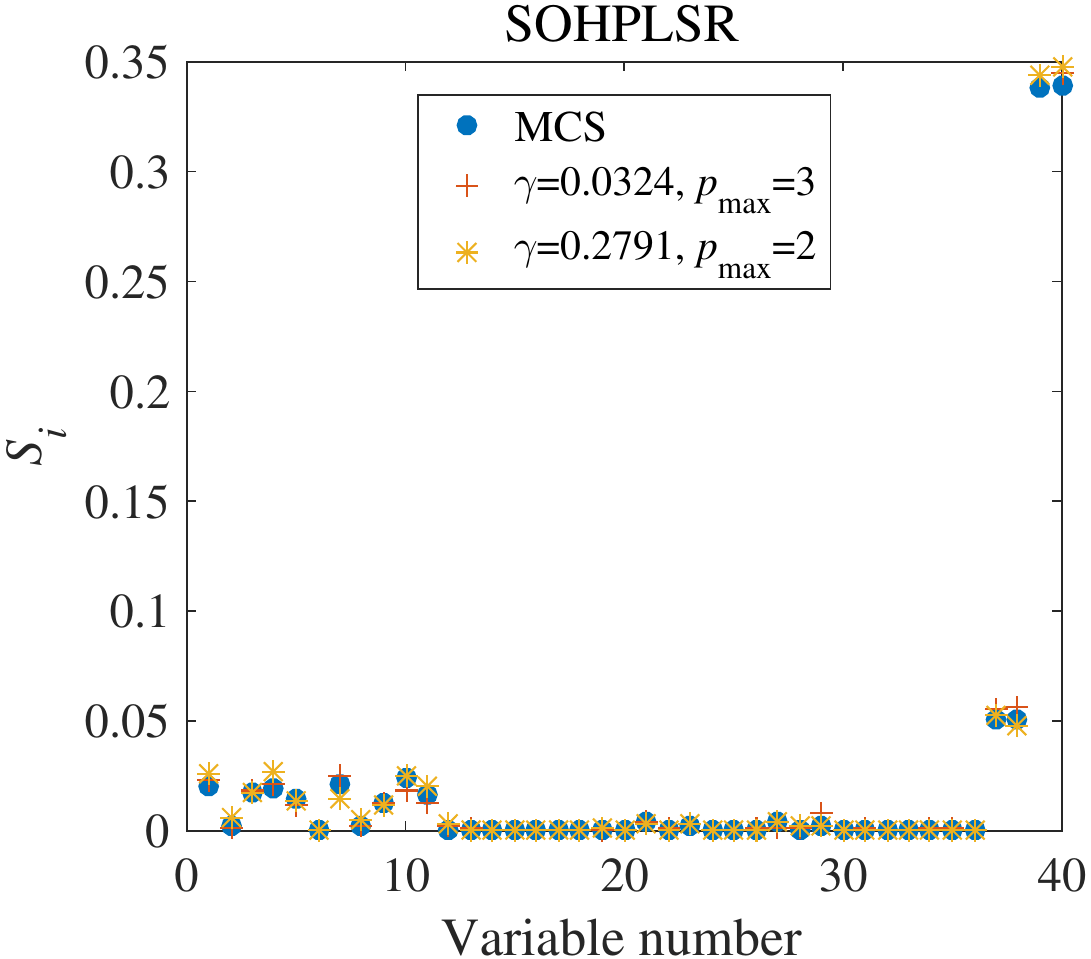}}}
	\caption{Comparison of the main Sobol indices}\label{Ex6_SiOLSHPLS}
	\vspace{0em}
\end{figure}

\begin{figure}[htbp]
	\centering
	\noindent\makebox[\textwidth][c]{	
		\subfigure[OLSR-PCE]
		{\includegraphics[width=0.5\textwidth]{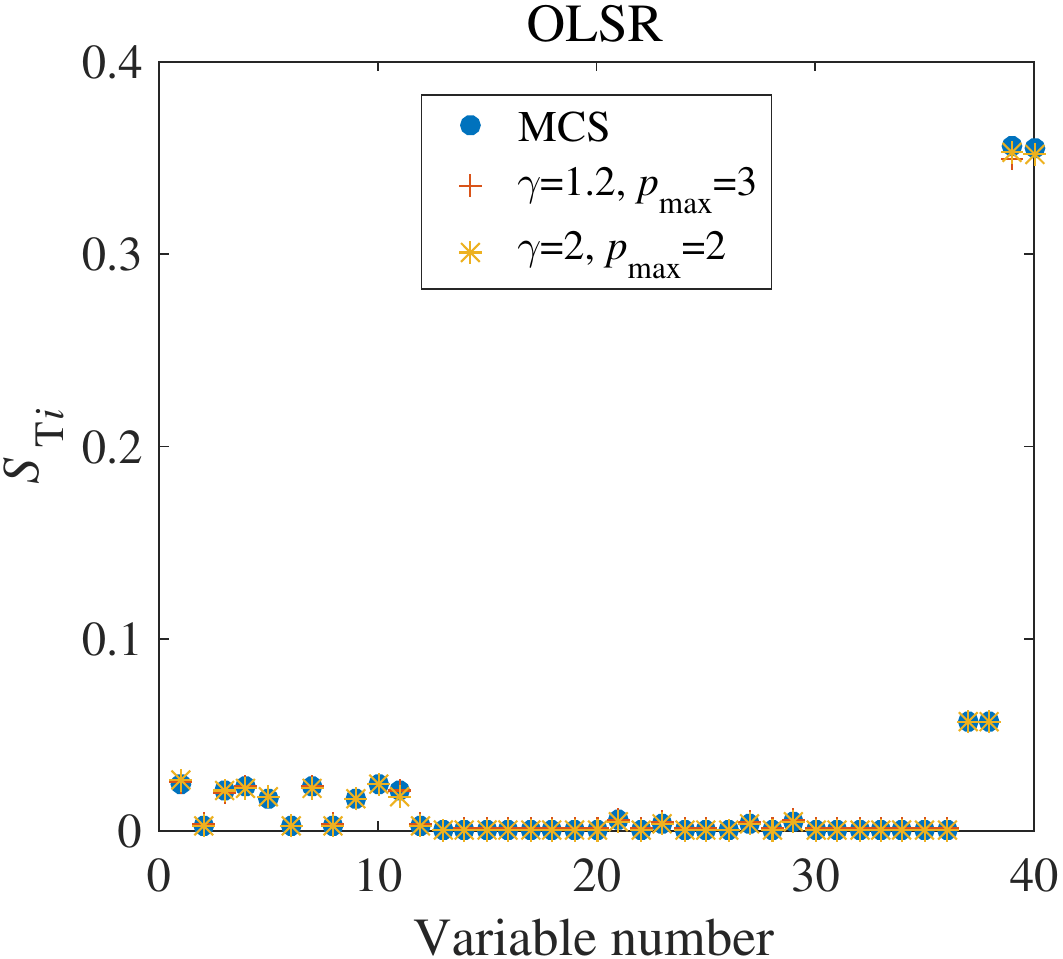}}
		\subfigure[SOHPLSR-PCE]
		{\includegraphics[width=0.5\textwidth]{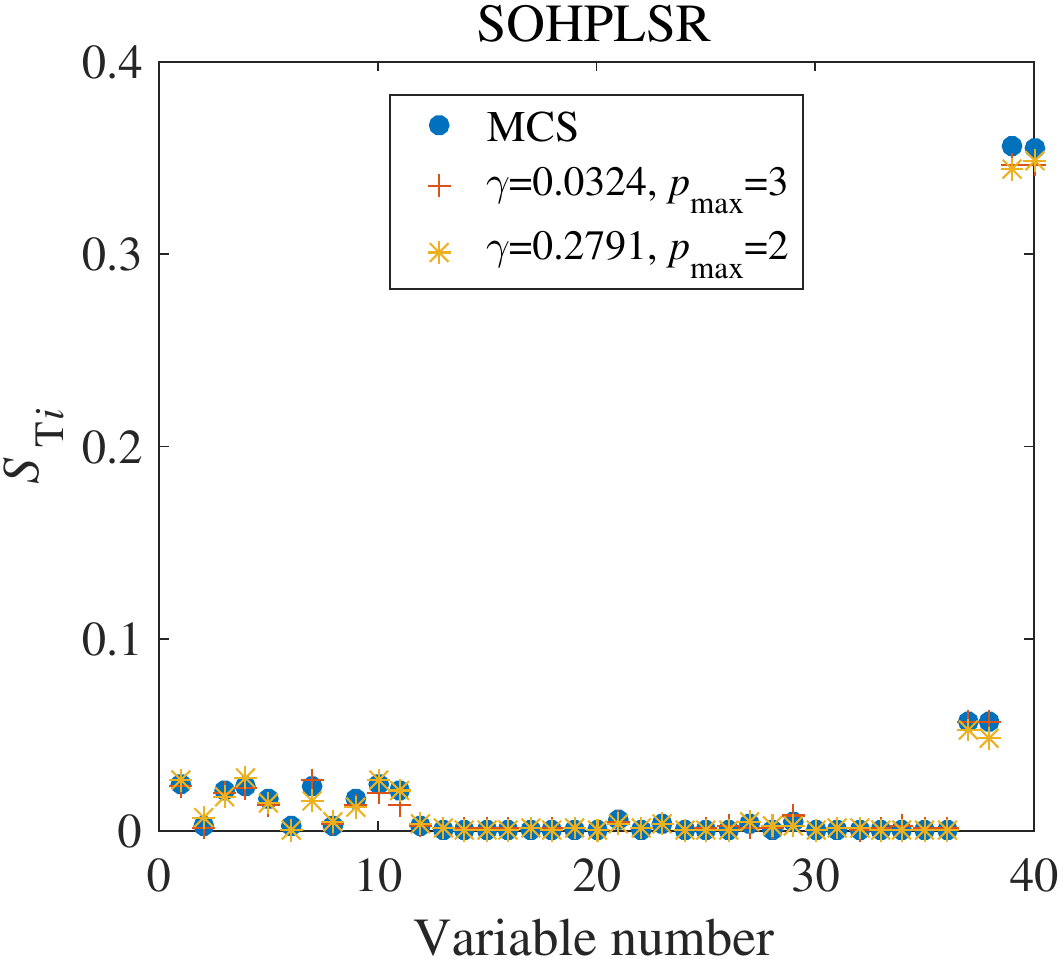}}}
	\caption{Comparison of the total Sobol indices}\label{Ex6_STOLSHPLS}
	\vspace{0em}
\end{figure}

The reference solution of each Sobol index is obtained by using MCS with $2\times 10^6$ samples. Comparisons of the main and total Sobol indices are illustrated in Figures \ref{Ex6_SiOLSHPLS} and \ref{Ex6_STOLSHPLS}.
It can be seen that variables 1, 3, 4, 5, 7, 9, 10, 11, 37-40 individually have significant impact on the uncertainty of model output, but their interactions have a weak impact. For OLSR, considerably accurate results of both $S_i$ and $S_{\mathrm{T}i}$ values can be obtained under $p_{\mathrm{max}}=2$ while the computational cost is only $1720/14808= 11.62\%$ as much as that of $p_{\mathrm{max}}=3$. For SOHPLSR, satisfactory results of both $S_i$ and $S_{\mathrm{T}i}$ values can be obtained with $N=240$ which is also smaller than that of $p_{\mathrm{max}}=3\ (N=400)$. In this example, third-order terms contribute little to the variance of the model output. However, to identify their contribution with acceptable accuracy, more samples are needed when $p_{\mathrm{max}}=3$. The computational efficiencies of SOPHLSR compared to that of OLSR are $1720/240= 7.1667$ for $p_{\mathrm{max}}=2$ and $14808/400= 37.02$ for $p_{\mathrm{max}}=3$, respectively.

\textit{Step 3}: Reliability analysis

\begin{figure}[htbp]
	\centering
	\includegraphics[width = 0.5\textwidth]{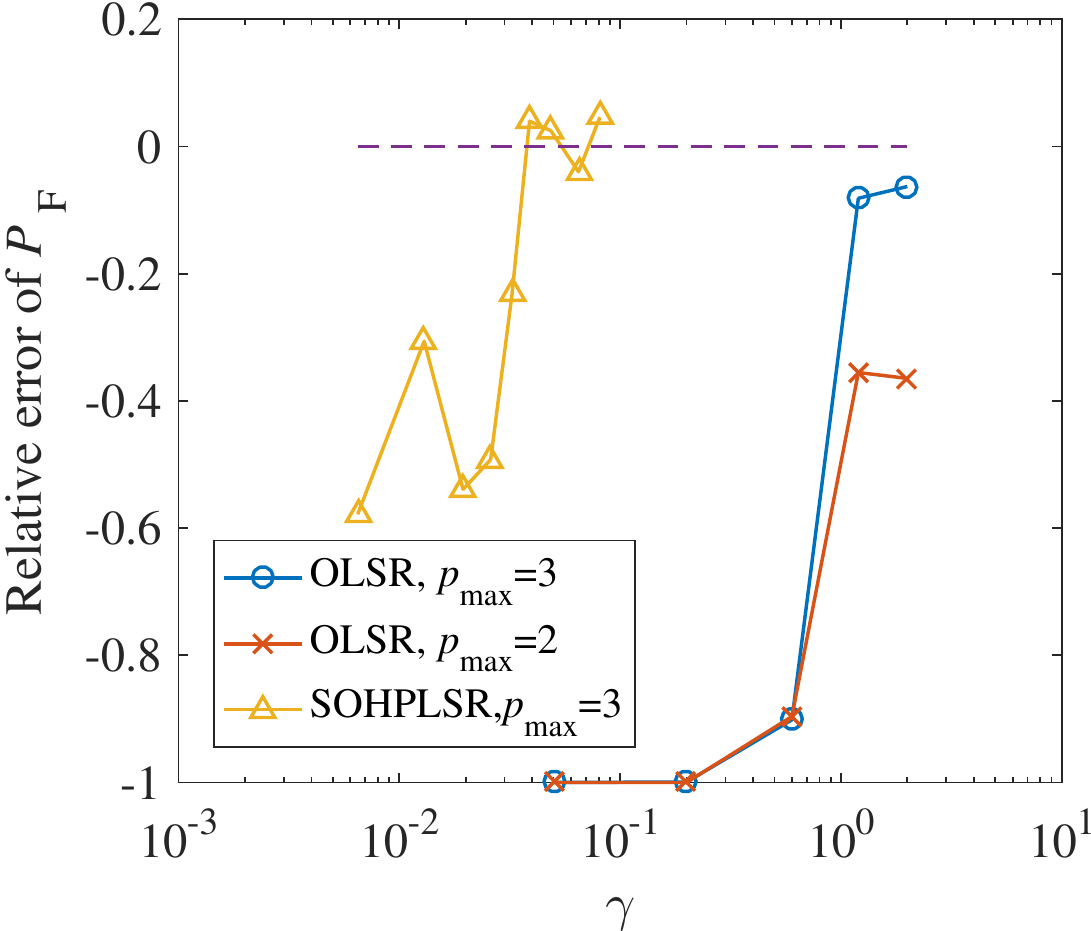}
	\caption{Comparison of relative errors of the failure probability}\label{Ex6_ePF_HPLS_AMM}
\end{figure}

The reference result of failure probability is 0.0051, which is obtained by using MCS with $3\times 10^6$ samples. The threshold for screening the important random inputs is set as 0.01. Relative errors of the failure probability computed by using OLSR under $p_{\mathrm{max}}=2$ and $p_{\mathrm{max}}=3$ are illustrated in Figure \ref{Ex6_ePF_HPLS_AMM}. 

It can been that the failure probabilities computed by OLSR under $p_{\mathrm{max}}=2$ fails to converge to the reference result since tail probability distribution of the model output is more sensitive to higher order statistics. Thus, to accurately compute the failure probability, higher values of $p_{\mathrm{max}}$ and $\gamma$ are needed. For SOHPLSR under $p_{\mathrm{max}}=3$, three more groups of experiment ($\phi$=15, 20, 25;\ $N$= 600, 800, 1000;\ $\gamma$= 0.0486, 0.0648, 0.0810) are added. The results are also illustrated in Figure \ref{Ex6_ePF_HPLS_AMM}. It can be seen that the failure probabilities computed by the proposed method converge to the reference result when $N$ reaches 600, indicating the computational gain is $24680/600\approx 41.1333$ compared to the OLSR-based counterpart.

\section{Conclusions}
This paper develops a novel non-intrusive method called SOHPLSR-PCE for global sensitivity and reliability analyses of high-dimensional models. The SOHPLSR-PCE shed new light on PCE by using hierarchical modeling and latent variable extraction. A state-of-the-art regression technique named PLSR is introduced to compute the latent factors that capture the most probabilistic information of polynomials at different variable levels. From the perspective of HDMR, the SOHPLSR-PCE automatically estimates the optimal interaction degree and the corresponding nonlinearity degrees. Three finite element models with different structural types and effective stochastic dimensions are employed to compare the relative performance of the prosed method and the traditional counterpart. The results demonstrate that the proposed method has the following properties: (1) For models with low effective stochastic dimensions (e.g. 2), the computational efficiencies of both the global sensitivity indices and the failure probabilities can be considerably high (e.g. 2 orders of magnitude higher than the traditional counterpart) without additional computational cost; (2) For models with moderate effective stochastic dimensions (e.g. 12-15), if only global sensitivity indices are needed, the computational efficiency is still rather high (e.g. tens of times higher than the traditional counterpart). When failure probabilities are needed, although more samples are required to keep accuracy, the computational efficiency is still much higher (e.g. ten times) than that of the traditional counterpart. (3) For the selection of $p_{\mathrm{max}}$, it is feasible to select a low value (e.g. 2) if only global sensitivity indices are needed. A higher value (e.g. 3) should be selected when reliability is to be analyzed. In summary, the proposed method not only has the potential of improving the computational efficiencies of global sensitivity and reliability analyses, but is promising for uncovering the latent hierarchical low dimensional structure of the models as well.

The proposed method has been so far applied to the models with univariate output and weak nonlinearities. Advanced methods for models with multivariate outputs and high nonlinearities are worth future investigations.

\section*{Acknowledgments}
This research was supported by the National Natural Science Foundation of China (NSFC, 51308158), and the China Postdoctoral Science Foundation (CPSF, 2013M541390), which are gratefully acknowledged by the authors. The valuable suggestions provided by Professor Guang-Chun Zhou are also gratefully acknowledged.
\section*{References}
\bibliography{SOPLSPCE}

\begin{thebibliography}{10}
\expandafter\ifx\csname url\endcsname\relax
  \def\url#1{\texttt{#1}}\fi
\expandafter\ifx\csname urlprefix\endcsname\relax\def\urlprefix{URL }\fi
\expandafter\ifx\csname href\endcsname\relax
  \def\href#1#2{#2} \def\path#1{#1}\fi

\bibitem{Blatman2010a}
G.~Blatman, B.~Sudret, Efficient computation of global sensitivity indices
  using sparse polynomial chaos expansions, Reliability Engineering \& System
  Safety 95~(11) (2010) 1216--1229.

\bibitem{Saltelli2010}
A.~Saltelli, P.~Annoni, I.~Azzini, F.~Campolongo, M.~Ratto, S.~Tarantola,
  Variance based sensitivity analysis of model output. design and estimator for
  the total sensitivity index, Computer Physics Communications 181~(2) (2010)
  259--270.

\bibitem{Zhang2014}
X.~Zhang, M.~D. Pandey, An effective approximation for variance-based global
  sensitivity analysis, Reliability Engineering \& System Safety 121 (2014)
  164--174.

\bibitem{Borgonovo2007}
E.~Borgonovo, A new uncertainty importance measure, Reliability Engineering \&
  System Safety 92~(6) (2007) 771--784.

\bibitem{Greegar2015}
G.~Greegar, C.~S. Manohar, Global response sensitivity analysis using
  probability distance measures and generalization of sobol's analysis,
  Probabilistic Engineering Mechanics 41 (2015) 21--33.

\bibitem{Sobol2009}
I.~M. Sobol, S.~Kucherenko, Derivative-based global sensitivity measures and
  their link with global sensitivity indices, Mathematics \& Computers in
  Simulation 79~(10) (2009) 3009--3017.

\bibitem{Sudret2015}
B.~Sudret, C.~V. Mai, Computing derivative-based global sensitivity measures
  using polynomial chaos expansions, Reliability Engineering \& System Safety
  134 (2015) 241--250.

\bibitem{Au1999}
S.~K. Au, J.~L. Beck, A new adaptive importance sampling scheme for reliability
  calculations, Structural safety 21~(2) (1999) 135--158.

\bibitem{Dai2012}
H.~Dai, H.~Zhang, W.~Wang, A support vector density-based importance sampling
  for reliability assessment, Reliability Engineering \& System Safety 106
  (2012) 86--93.

\bibitem{Dai2016}
H.~Dai, H.~Zhang, W.~Wang, A new maximum entropy-based importance sampling for
  reliability analysis, Structural Safety 63 (2016) 71--80.

\bibitem{Papaioannou2015}
I.~Papaioannou, W.~Betz, K.~Zwirglmaier, D.~Straub, Mcmc algorithms for subset
  simulation, Probabilistic Engineering Mechanics 41 (2015) 89--103.

\bibitem{Zuev2015}
K.~Zuev, Subset simulation method for rare event estimation: An introduction,
  arXiv preprint arXiv:1505.03506.

\bibitem{Koutsourelakis2004}
P.~S. Koutsourelakis, Reliability of structures in high dimensions. part ii.
  theoretical validation, Probabilistic engineering mechanics 19~(4) (2004)
  419--423.

\bibitem{Koutsourelakis2004a}
P.~S. Koutsourelakis, H.~J. Pradlwarter, G.~I. Schueller, Reliability of
  structures in high dimensions, part i: algorithms and applications,
  Probabilistic Engineering Mechanics 19~(4) (2004) 409--417.

\bibitem{Konakli2016}
K.~Konakli, B.~Sudret, Reliability analysis of high-dimensional models using
  low-rank tensor approximations, Probabilistic Engineering Mechanics 46 (2016)
  18--36.

\bibitem{Dai2015}
H.~Dai, H.~Zhang, W.~Wang, A multiwavelet neural network‐based response
  surface method for structural reliability analysis, Computer aided Civil and
  Infrastructure Engineering 30~(2) (2015) 151--162.

\bibitem{Dai2017}
H.~Dai, Z.~Cao, A wavelet support vector machine-based neural network metamodel
  for structural reliability assessment, Computer-Aided Civil and
  Infrastructure Engineering 32~(4) (2017) 344--357.

\bibitem{Wiener1938}
N.~Wiener, The homogeneous chaos, American Journal of Mathematics 60~(4) (1938)
  897--936.

\bibitem{Ghanem2003}
R.~G. Ghanem, P.~D. Spanos, Stochastic finite elements: a spectral approach,
  Courier Corporation, 2003.

\bibitem{Blatman2009}
G.~Blatman, Adaptive sparse polynomial chaos expansions for uncertainty
  propagation and sensitivity analysis, Ph.D. thesis, Clermont-Ferrand 2
  (2009).

\bibitem{Blatman2010}
G.~Blatman, B.~Sudret, An adaptive algorithm to build up sparse polynomial
  chaos expansions for stochastic finite element analysis, Probabilistic
  Engineering Mechanics 25~(2) (2010) 183--197.

\bibitem{Abraham2017}
S.~Abraham, M.~Raisee, G.~Ghorbaniasl, F.~Contino, C.~Lacor, A robust and
  efficient stepwise regression method for building sparse polynomial chaos
  expansions, Journal of Computational Physics 332 (2017) 461--474.

\bibitem{Liu2018}
Q.~Liu, X.~Zhang, X.~Huang, A sparse surrogate model for structural reliability
  analysis based on the generalized polynomial chaos expansion, Proceedings of
  the Institution of Mechanical Engineers, Part O: Journal of Risk and
  Reliability\href {http://dx.doi.org/10.1177/1748006X18804047}
  {\path{doi:10.1177/1748006X18804047}}.

\bibitem{Blatman2011a}
G.~Blatman, B.~Sudret, Adaptive sparse polynomial chaos expansion based on
  least angle regression, Journal of Computational Physics 230~(6) (2011)
  2345--2367.

\bibitem{Cheng2018}
K.~Cheng, Z.~Lu, Adaptive sparse polynomial chaos expansions for global
  sensitivity analysis based on support vector regression, Computers \&
  Structures 194 (2018) 86--96.

\bibitem{Cheng2018a}
K.~Cheng, Z.~Lu, Sparse polynomial chaos expansion based on d-morph regression,
  Applied Mathematics and Computation 323 (2018) 17--30.

\bibitem{Jakeman2017}
J.~D. Jakeman, A.~Narayan, T.~Zhou, A generalized sampling and preconditioning
  scheme for sparse approximation of polynomial chaos expansions, SIAM Journal
  on Scientific Computing 39~(3) (2017) A1114--A1144.

\bibitem{Rosipal2006}
R.~Rosipal, N.~Kr{\"a}mer, Overview and recent advances in partial least
  squares, in: Subspace, Latent Structure and Feature Selection, Springer,
  2006, pp. 34--51.

\bibitem{Rosipal2010}
R.~Rosipal, Nonlinear partial least squares: An overview, Chemoinformatics and
  advanced machine learning perspectives: complex computational methods and
  collaborative techniques (2010) 169--189.

\bibitem{Zhao2019}
W.~Zhao, L.~Bu, Global sensitivity analysis with a hierarchical sparse
  metamodeling method, Mechanical Systems and Signal Processing 115 (2019)
  769--781.

\bibitem{Soize2004}
C.~Soize, R.~Ghanem, Physical systems with random uncertainties: chaos
  representations with arbitrary probability measure, SIAM Journal on
  Scientific Computing 26~(2) (2004) 395--410.

\bibitem{Rahman2017}
S.~Rahman, Wiener–-hermite polynomial expansion for multivariate gaussian
  probability measures, Journal of Mathematical Analysis and Applications
  454~(1) (2017) 303--334.

\bibitem{Cameron1947}
R.~H. Cameron, W.~T. Martin, The orthogonal development of non-linear
  functionals in series of fourier-hermite functionals, Annals of Mathematics
  48~(48) (1947) 385--392.

\bibitem{Smolyak1963}
S.~Smolyak, Quadrature and interpolation formulas for tensor products of
  certain classes of functions, in: Soviet Math. Dokl., Vol.~4, 1963, pp.
  240--243.

\bibitem{Mehmood2012}
T.~Mehmood, K.~H. Liland, L.~Snipen, S.~S{\ae}b{\o}, A review of variable
  selection methods in partial least squares regression, Chemometrics \&
  Intelligent Laboratory Systems 118~(16) (2012) 62--69.

\bibitem{Cao2008}
K.-A. L{\^e}~Cao, D.~Rossouw, C.~Robert-Grani{\'e}, P.~Besse, A sparse pls for
  variable selection when integrating omics data, Statistical applications in
  genetics and molecular biology 7~(1) (2008) Article 35.

\bibitem{Rabitz1999b}
H.~Rabitz, {\"O}.~F. Ali{\c{s}}, General foundations of high-dimensional model
  representations, Journal of Mathematical Chemistry 25~(2-3) (1999) 197--233.

\bibitem{Sudret2008}
B.~Sudret, Global sensitivity analysis using polynomial chaos expansions,
  Reliability Engineering \& System Safety 93~(7) (2008) 964--979.

\bibitem{Rahman2008}
S.~Rahman, A polynomial dimensional decomposition for stochastic computing,
  International Journal for Numerical Methods in Engineering 76~(13) (2008)
  2091--2116.

\bibitem{Kattan2008}
P.~I. Kattan, MATLAB Guide to Finite Elements, Springer Berlin Heidelberg,
  2008.

\end{thebibliography}
\end{document}